\documentclass{gtart_a}
\pdfoutput=1

\usepackage{graphicx}
\usepackage[percent]{overpic}
\usepackage{placeins}
\usepackage{units}

\title{Criticality for the Gehring link problem}

\author[Cantarella]{Jason Cantarella}
\givenname{Jason}
\surname{Cantarella}
\address{[JC, JHGF]\ \ Department of Mathematics\\
University of Georgia\\\newline
Athens, GA 30602\\USA}
\email{cantarel@math.uga.edu}
\urladdr{}

\author[Fu]{Joseph H\,G Fu}
\givenname{Joseph H\,G}
\surname{Fu}
\email{fu@math.uga.edu}
\urladdr{}

\author[Kusner]{Rob Kusner}
\givenname{Rob}
\surname{Kusner}
\address{Department of Mathematics\\
University of Massachusetts\\\newline
Amherst, MA 01003\\USA}
\email{kusner@math.umass.edu}
\urladdr{}

\author[Sullivan]{John M Sullivan}
\givenname{John}
\surname{Sullivan}
\address{Institut fuer Mathematik\\
Technische Universitat Berlin\\\newline
DE--10623 Berlin\\Germany}
\email{sullivan@math.tu-berlin.de}
\urladdr{}

\author[Wrinkle]{Nancy C Wrinkle}
\givenname{Nancy}
\surname{Wrinkle}
\address{Department of Mathematics\\
Northeastern Illinois University\\\newline
Chicago, IL 60625\\USA}
\email{N-Wrinkle@neiu.edu}
\urladdr{}

\volumenumber{10}
\issuenumber{}
\publicationyear{2006}
\papernumber{45}
\lognumber{0610}
\startpage{2055}
\endpage{2115}

\doi{}
\MR{}
\Zbl{}

\arxivreference{math.DG/0402212}


\received{16 May 2005}
\revised{10 October 2006}
\accepted{23 May 2006}
\published{14 November 2006}
\publishedonline{14 November 2006}
\proposed{Yasha Eliashberg}
\seconded{Joan Birman, Tobias Colding}
\corresponding{}
\editor{Michael Freedman,Joan Birman}
\version{}

\AtBeginDocument{\let\bar\wbar\let\tilde\wtilde\let\hat\what}
\AtBeginDocument{\renewcommand{\phi}{\varphi}}
\AtBeginDocument{\renewcommand{\S}{\mathbb{S}}}
\AtBeginDocument{\renewcommand{\H}{\mathcal{H}}}
\AtBeginDocument{\renewcommand{\d}{\mathrm{d}}}

\let\tsty\textstyle

\newcommand{\nstartwo}[1]{$#1$$*$$2$}
\newcommand{\arc}[1]{#1} 

\newcommand{\approaches}{\to}
\newcommand{\cross}{\times}

\newcommand{\disj}{\sqcup} 
\newcommand{\bigdisj}{\bigsqcup}

\newcommand{\dist}{\mathrm{Dist}}

\newcommand{\len}{\mathrm{Len}}
\newcommand{\eps}{\varepsilon}

\newcommand{\isom}{\cong}

\newcommand{\bd}{\partial}

\newcommand{\K}{\mathcal{K}}

\newcommand{\UPL}{L^{(2)}}
\newcommand{\Aw}{A_\mathrm{W}}
\newcommand{\Ag}{A_\mathrm{S}}

\newcommand{\ds}{\,\d s}
\newcommand{\dx}{\,\d x}
\newcommand{\dz}{\,\d z}
\newcommand{\du}{\,\d u}
\newcommand{\dmu}{\,\d\mu}

\newcommand{\dplus}{\delta^+} 

\newcommand{\ddtz}[3]{\frac{\d #1}{\d t^{#2}} {#3} \Big|_{t=0}} 
\newcommand{\ddtzp}[1]{\ddtz{}{}{#1}} 

\newcommand{\tild}[1]{\tilde{#1}}
\newcommand{\ttild}[1]{\hat{#1}}
\newcommand{\tg}{{\tild\gamma}}

\newcommand{\FTC}{\operatorname{FTC}}

\newcommand{\gThi}{\operatorname{LThi}} 
\newcommand{\gStrut}{\operatorname{Strut}}
\newcommand{\Gehring}{Link-}
\newcommand{\gehring}{link-}

\newcommand{\mass}{\operatorname{mass}}

\newcommand{\VF}{\operatorname{VF}}
\newcommand{\wall}{\operatorname{Wall}}
\newcommand\grad{\nabla}
\newcommand{\Len}{\operatorname{Len}}

\makeatletter
\def\cnewtheorem#1[#2]#3{\newtheorem{#1}{#3}[section]
\expandafter\let\csname c@#1\endcsname\c@theorem}

\newtheorem{theorem}{Theorem}[section]
\cnewtheorem{lemma}[theorem]{Lemma}
\cnewtheorem{corollary}[theorem]{Corollary}
\cnewtheorem{proposition}[theorem]{Proposition}

\theoremstyle{definition}
\cnewtheorem{example}[theorem]{Example}
\cnewtheorem{remark}[theorem]{Remark}
\newtheorem*{definition}{Definition}

\newtheorem*{question}{Question}

\makeatother  

\numberwithin{equation}{section}

\newcommand{\lem}[1]{\fullref{lem:#1}}
\newcommand{\thm}[1]{\fullref{thm:#1}}
\newcommand{\prop}[1]{\fullref{prop:#1}}
\newcommand{\cor}[1]{\fullref{cor:#1}}
\newcommand{\xmpl}[1]{\fullref{ex:#1}}
\newcommand{\eqn}[1]{\eqref{eq:#1}}

\def\figr#1{\fullref{fig:#1}}
\def\secn#1{\fullref{sec:#1}}

\newcommand{\incgr}[2]{\includegraphics[scale=#1]{#2}}

\newcommand{\figcap}[3]{\caption[#2]{#3}\label{fig:#1}}

\newcommand{\figs}[5]{%
  \begin{figure}[ht!]
    \centerline{\incgr{#2}{#3}}
    \figcap{#1}{#4}{#5}
  \end{figure}
}

\newcommand{\figthree}[7]{%
  \begin{figure}[ht!]
    \incgr{#2}{#3} \hfill \incgr{#2}{#4} \hfill \incgr{#2}{#5}
    \figcap{#1}{#6}{#7}
  \end{figure}
}

\renewcommand{\theenumi}{\alph{enumi}}

\newenvironment{enuma} {\begin{enumerate}}{\end{enumerate}}
\newcommand{\itm}[1]{{\eqref{itm:#1}}}

\begin{document}

\begin{asciiabstract}
In 1974, Gehring posed the problem of minimizing the length of two
linked curves separated by unit distance.  This constraint can be
viewed as a measure of thickness for links, and the ratio of length
over thickness as the ropelength.  In this paper we refine
Gehring's problem to deal with links in a fixed link-homotopy class:
we prove ropelength minimizers exist and introduce a theory of
ropelength criticality.

Our balance criterion is a set of necessary and sufficient conditions
for criticality, based on a strengthened, infinite-dimensional version
of the Kuhn--Tucker theorem.  We use this to prove that every critical
link is C^1 with finite total curvature.  The balance criterion also
allows us to explicitly describe critical configurations (and presumed
minimizers) for many links including the Borromean rings.  We also
exhibit a surprising critical configuration for two clasped ropes:
near their tips the curvature is unbounded and a small gap appears between
the two components.  These examples reveal the depth
and richness hidden in Gehring's problem and our natural
extension.
\end{asciiabstract}

\begin{webabstract}
In 1974, Gehring posed the problem of minimizing the length of two
linked curves separated by unit distance.  This constraint can be
viewed as a measure of thickness for links, and the ratio of length
over thickness as the ropelength.  In this paper we refine
Gehring's problem to deal with links in a fixed link-homotopy class:
we prove ropelength minimizers exist and introduce a theory of
ropelength criticality.

Our balance criterion is a set of necessary and sufficient conditions
for criticality, based on a strengthened, infinite-dimensional version
of the Kuhn--Tucker theorem.  We use this to prove that every critical
link is $C^1$ with finite total curvature.  The balance criterion also
allows us to explicitly describe critical configurations (and presumed
minimizers) for many links including the Borromean rings.  We also
exhibit a surprising critical configuration for two clasped ropes:
near their tips the curvature is unbounded and a small gap appears between
the two components.  These examples reveal the depth
and richness hidden in Gehring's problem and our natural
extension.
\end{webabstract}

\begin{htmlabstract}
<p class="noindent">
In 1974, Gehring posed the problem of minimizing the length of two
linked curves separated by unit distance.  This constraint can be
viewed as a measure of thickness for links, and the ratio of length
over thickness as the ropelength.  In this paper we refine
Gehring's problem to deal with links in a fixed link-homotopy class:
we prove ropelength minimizers exist and introduce a theory of
ropelength criticality.
</p>
<p class="noindent">
Our balance criterion is a set of necessary and sufficient conditions
for criticality, based on a strengthened, infinite-dimensional version
of the Kuhn&ndash;Tucker theorem.  We use this to prove that every critical
link is&nbsp;C<sup>1</sup> with finite total curvature.  The balance criterion also
allows us to explicitly describe critical configurations (and presumed
minimizers) for many links including the Borromean rings.  We also
exhibit a surprising critical configuration for two clasped ropes:
near their tips the curvature is unbounded and a small gap appears between
the two components.  These examples reveal the depth
and richness hidden in Gehring's problem and our natural
extension.
</p>
\end{htmlabstract}

\begin{abstract}
In 1974, Gehring posed the problem of minimizing the length of two
linked curves separated by unit distance.  This constraint can be
viewed as a measure of thickness for links, and the ratio of length
over thickness as the ropelength.  In this paper we refine
Gehring's problem to deal with links in a fixed link-homotopy class:
we prove ropelength minimizers exist and introduce a theory of
ropelength criticality.

Our balance criterion is a set of necessary and sufficient conditions
for criticality, based on a strengthened, infinite-dimensional version
(\thm{ourkt})
of the Kuhn--Tucker theorem.  We use this to prove that every critical
link is~$C^1$ with finite total curvature.  The balance criterion also
allows us to explicitly describe critical configurations (and presumed
minimizers) for many links including the Borromean rings.  We also
exhibit a surprising critical configuration for two clasped ropes:
near their tips the curvature is unbounded and a small gap appears between
the two components.  These examples reveal the depth
and richness hidden in Gehring's problem and our natural
extension.
\end{abstract}

\maketitle

\section{Introduction}
\begin{flushright}
\emph{Suppose that~$A$ and~$B$ are disjoint linked Jordan curves in~$\R^3$
\\ which lie at a distance~$1$ from each other.
\\ Show that the length of~$A$ is at least~$2\pi$.}

---Fred Gehring, 1974
\end{flushright}
Gehring's problem, which appeared in a conference proceedings~\cite{MR52:5944},
was soon solved by Marvin Ortel. Because Ortel's elegant solution was never
published, we reproduce it here with his permission:  
Fix any point $a \in A$; the cone on~$A$ from~$a$ is a disk spanning~$A$.
Since~$A$ and~$B$ are linked, $B$ meets this disk at some point $b\in B$,
lying on a chord of~$A$.  Because $\dist(A,b) \ge 1$, projecting~$A$ 
to the unit sphere~$S$ around~$b$ does not increase its length. The
projection is a closed curve joining two antipodal points on~$S$,
and so has length at least~$2\pi$.
(Further proofs and generalizations to linked spheres in higher
dimensions were published by Edelstein and Schwarz \cite{EdSc}, Osserman \cite{Oss} and Gage \cite{MR82b:52017,MR82j:52024}.)

The unique minimizing configuration for Gehring's problem is a Hopf link
consisting of two congruent circles in perpendicular planes, each passing
through the other's center.
This leads to a natural question: what are the length-minimizing shapes of
other link types when the different components stay unit distance apart?
This constraint prevents different components from crossing
each other, but we cannot expect to fix the link type exactly.
Instead, the natural setting for this problem is Milnor's notion
of link homotopy: two links are link-homotopic if one can be deformed
into the other while keeping different components disjoint.
Clearly one link can be deformed into another while keeping
all components at unit distance if and only if they are link homotopic.

We will define the \emph{\gehring thickness\/} of a link to be the
minimum distance between different components.  The problem
we consider is then to minimize length in a link-homotopy class,
subject to the constraint of fixed \gehring thickness.
Equivalently, we could minimize the \emph{\gehring ropelength\/}
of the link, meaning the quotient of length over thickness.

In~\cite{CKS2}, we found length-minimizing links
under a similar constraint: that a normal tube of diameter one
around the link stay embedded.  It is easy to see that the examples
constructed there (like the one in \figr{known})
are also global minima (in their respective link-homotopy classes)
for the Gehring problem.  The focus of
this paper will be on \emph{critical\/} configurations.
Our main result is a balance criterion (\thm{gbalance}, \mbox{\cor{gehr-bal}}),
which states that a link is \gehring ropelength critical
if and only if the tension force in the curve is balanced
by a system of compressive forces between pairs of points on
different components of~$L$ realizing the minimum distance.

This balance criterion is based mainly on an improved, infinite-dimensional
version (\thm{ourkt}) of the Kuhn--Tucker theorem on constrained optimization,
which is essentially a very general method of Lagrange multipliers.
The other key technical element is a careful application
of Clarke's differentiation theorem for min-functions
(\thm{clarke}).

The direct method
shows that there is a (rectifiable) minimizer for \gehring ropelength in
each link-homotopy class.  An interesting problem is to determine the
regularity of these minimizers or other critical points.  The previously
known minimizers were~$C^{1,1}$ but not~$C^2$.  Our balance criterion
allows us to prove that all \gehring ropelength-critical curves are~$C^1$ with
finite total curvature (\fullref{prop:gcrit-c1}).

We next consider generalized links, which may include open components with
constrained endpoints, or which avoid fixed obstacles.  After extending
our balance criterion and existence results to this setting,
we analyze the problem of the \emph{simple clasp\/}.  A clasp consists
of two linked arcs whose endpoints are constrained to parallel planes
(as in \figr{clasp}).
A generalization to clasps of different opening angles provides a model
for the strands of rope in a woven cloth or net.
The balance criterion lets us construct explicit critical configurations
(\thm{clasp}) of these generalized links; we conjecture
they are the length-minimizers subject to the constraint that
the arcs remain at unit distance from each other.
Our critical clasp has a number of surprising features,
including a point of infinite curvature and a small gap (at the center
of the clasp) between the tubes around the two components.
This configuration is $C^{1,\nicefrac{2}{3}}$ and may represent the
worst regularity of any critical curve.

We end by constructing a ropelength-critical configuration
(and presumed minimizer) for the Borromean rings.
In all the other known critical configurations for closed links,
each component is a convex plane curve built from straight
segments and arcs of circles.  In our Borromean rings,
the components are still planar, but are nonconvex and are built
from different pieces including parts of a clasp curve.  In a sense,
this is the first nontrivial example of a ropelength-critical link.

Our methods will have a number of other applications. In particular, we
have used them to describe critical configurations for the ``standard''
ropelength problem for knots and links: minimize the length of a~$C^1$ link
subject to the constraint that the normal neighborhood of unit diameter
remains embedded.  We will publish these results in a sequel~\cite{CFKSW2}
to the current paper.  We can also consider minimization not of length
but of other objective functions like elastic bending energy, again
subject to a thickness constraint.
Analogs of our balance criterion may be useful in
describing other flexible mechanisms, such as thick surfaces.

We note that von der Mosel and Schuricht~\cite{vdms} have used a similar
approach (via Clarke's theorem and a functional-analytic version of
Lagrange multipliers) to derive necessary, but not sufficient, conditions
for criticality for the ropelength functional of~\cite{CKS2}. We will treat
the same functional in our forthcoming sequel~\cite{CFKSW2} and will
offer a comparison of the two methods there.  We also note that
Starostin has given~\cite{starostin-forma} an independent derivation 
of the tight clasp of \fullref{sec:gehr-clasp},
though he does not prove that it is critical.

\section{\Gehring thickness for closed links} \label{sec:gehringdef}

In order to reformulate Gehring's problem,
we first establish some basic terminology.
Remember that a compact, oriented $1$--manifold-with-boundary~$M$
is a finite union of components,
each of which is homeomorphic to a circle~$\S^1$ or an interval~$[0,1]$.

\begin{definition}
A \emph{parametrized curve\/} is a mapping from
a compact, oriented $1$--manifold-with-boundary~$M$ to~$\R^3$.
Two parametrized curves are equivalent if they differ by
an orientation-preserving reparametrization (thatis, by composition
with an orientation-preserving self-homeomorphism of~$M$).
A \emph{curve\/}~$L$ in~$\R^3$ is an equivalence class of parametrized curves.
We say~$L$ is \emph{closed\/} when
each component of its domain~$M$ is a circle,
that is, when its boundary~$\bd L$ is empty.
\end{definition}

Even though our curves may have self-intersections,
we will usually refer to points on the curve as if they were simply
points of its image in~$\R^3$. The meaning should be clear from context.

The \emph{length\/} $\len(L)$ of any curve~$L$ is defined to be the supremal
length of all polygons inscribed in~$L$. A curve has finite length,
or is \emph{rectifiable\/}, if and only if it has a Lipschitz (that is, $C^{0,1}$)
parametrization.  One such parametrization is then by arclength~$s$.
Any rectifiable curve has a well-defined unit tangent vector
$T=\d L/\d s$ almost everywhere.

\begin{definition} 
\label{defn:gthi}
The \emph{\gehring thickness\/} $\gThi(L)$ of a curve~$L$ is the minimum distance 
between points on different components of~$L$.
This is the supremal~$\eps$ for which the $(\nicefrac{\eps}{2})$--neighborhoods
of the components of~$L$ are disjoint.
\end{definition}

For now, we will consider only the case of closed curves,
where each component is a circle.  (We will deal with
generalized links---with endpoint constraints---later in \secn{genlink}.)
So suppose we start with a closed curve~$L$ and we want to minimize length
under the constraint that the \gehring thickness remains at least one.
Since we can rescale any link to have $\gThi \ge 1$,
this problem is the same as minimizing \emph{(\gehring) ropelength\/},
the quotient of length by \gehring thickness.

The thickness constraint naturally prevents different components
from passing through each other, but does not prevent any given
component from changing its knot type through self-intersections.
This is the setting for Milnor's work on link homotopy:

\begin{definition}
\label{defn:linkhomotopy}
A \emph{link\/} is a closed curve with disjoint components.
The \emph{link-homotopy\/} class of a link~$L$, denoted $[[L]]$, is the set of
curves homotopic to~$L$ through configurations that keep different
components of~$L$ disjoint.
\end{definition}

Note that, for our purposes, configurations of~$L$ where
some components have self-intersections are still considered to be links,
and are included in~$[[L]]$.

For two-component links, Milnor~\cite{MR17:70e} showed that linking
number is the only link-homotopy invariant.  For links of many components,
the topological situation is more complicated, but a complete
classification of links up to link homotopy was provided by Habegger and
Lin~\cite{MR91e:57015}.
We will prove in \secn{gehrbal} that in every link-homotopy class
there is a curve minimizing ropelength.
We show these minimizers are always~$C^1$, though our
examples suggest that they may not always have bounded curvature.

\section{The derivative of \gehring thickness} \label{sec:dgthi}
We want to define critical configurations of~$L$ subject to the
thickness constraint $\gThi(L) \ge 1$.
Because $\gThi$ is defined as the minimum of a collection
of distances between points on different components,
the equation $\gThi \ge 1$ acts like a collection of many constraints.
To make this notion precise, we will apply a theorem of Clarke
to compute the derivative of $\gThi$ as we vary the curve~$L$.

Given any curve~$L$, let~$\UPL$ be the compact
set of all unordered pairs $\{x,y\}$ of points
on distinct components of~$L$.  The \gehring thickness
of~$L$ is simply the minimum over~$\UPL$ of the distance
function $\dist\{x,y\} := |y-x|$.

We often want to consider a \emph{continuous
deformation}~$L_t$ of a curve~$L$: fixing any parametrization~$f$ of~$L$,
that means a continuous family~$f_t$ of parametrized curves with $f_0=f$.
(When we reparametrize~$L$, we apply the \emph{same\/} reparametrization
to~$L_t$ at all times~$t$.)  We assume that~$L_t$ is~$C^1$ in~$t$;
the initial velocity of~$L_t$ will then be
given by some (continuous, $\R^3$--valued) \emph{vector field\/}~$\xi$ along~$L$.
We let~$\VF(L)$ denote the space of all such vector fields.
Formally, these are sections of the bundle $f^*T\R^3$ pulled back from 
the tangent bundle of~$\R^3$ by the parametrization~$f$ of~$L$. 
Identifying any tangent space to~$\R^3$ with~$\R^3$ itself, this
is simply a map from the domain~$M$ to~$\R^3$.  Again, when we reparametrize
a curve~$L$, we apply the same reparametrization to any
vector field~$\xi$.

Consider a curve~$L$ with $\gThi(L)>0$.
If~$L_t$ is a continuous deformation of~$L$, with initial
velocity given by some $\xi\in\VF(L)$, then for each pair
$\{x,y\}\in \UPL$, we clearly have
$$\delta_\xi \dist\{x,y\} := \ddtzp{|y-x|}
 \!= \frac{\langle \xi_y - \xi_x, y - x \rangle}{|y-x|}.$$
(Even if~$L$ is not embedded, the condition $\gThi(L)>0$
implies~$x$ and~$y$ cannot coincide in~$\R^3$, so this
formula is always meaningful.)

A function like $\gThi$, defined as the minimum of a compact
family of smooth functions, is sometimes called a min-function.
Clarke's differentiation theorem for min-functions says
that---just as in the case when the compact family is finite---the
derivative of a min-function is the smallest derivative of those
smooth functions that achieve the minimum.
More precisely, specializing \cite[Theorem 2.1]{clarke}
to the case we need, we have:

\begin{theorem}[Clarke]
\label{thm:clarke}
Suppose for some compact space~$K$ and some $\eps>0$, we have a
family of~$C^1$ functions $f_k\co (-\eps,\eps)\to\R$, for $k\in K$.
Suppose further that~$f_k(t)$ and~$f'_k(t)$ are lower semicontinuous
on $K\cross(-\eps,\eps)$.  Let $f(t):=\min_{k\in K} f_k(t)$.
Then~$f$ has one-sided derivatives, and
$$\ddtz{f}{+}{} = \min_{k\in K_0} f'_k(0),$$
where $K_0 := \{k\in K: f_k(0)=f(0)\}$ is the set of~$k$ for which
the minimum in the definition of~$f$ is achieved when $t=0$.
\end{theorem}

To apply this theorem to thickness, suppose we have a variation~$L_t$
of the curve~$L$, and let $\xi\in\VF(L)$ be its initial velocity.
The \gehring thickness $\gThi(L)$ is written as a minimum over $K=\UPL$
of the pairwise distance.  Clarke's theorem picks out those pairs
achieving the minimum:~$K_0$ is the set of pairs achieving the minimum
distance $\gThi(L)$.

In rigidity theory, the vertices
of a tensegrity framework are joined by \emph{bars\/} whose
length is fixed, \emph{cables\/} whose length can shrink but not grow,
and \emph{struts\/} whose length can grow but not shrink (compare Roth and Whiteley \cite{MR82m:51018}).
Thus, we borrow the term ``strut'' to describe the pairs in~$K_0$:

\begin{definition}
\label{defn:gstrut}
An unordered pair of points $\{x,y\}$ on different components of~$L$ is a
\emph{strut\/} if $|y-x| = \gThi(L)$.
The space of all struts of~$L$ is denoted $\gStrut(L)\subset L^{(2)}$.
\end{definition}

Struts correspond to points of contact between tubes around
the different components of~$L$.
Our balance criterion will show how the segment~$\overline{xy}$
can be viewed as carrying a force pushing outward on its endpoints.

Applying Clarke's theorem to \gehring thickness, we get:
\begin{corollary}\label{cor:dgthi}
For any curve~$L$, and any variation vector field $\xi\in\VF(L)$,
the (one-sided) first variation of \gehring thickness is
$$\dplus_\xi \gThi(L) = \min_{\gStrut(L)}\!\! \delta_\xi \dist\{x,y\}.$$
\end{corollary}
Note that $\delta_\xi \dist\{x,y\}$ is a continuous function
of~$x$ and~$y$, and for any fixed $\{x,y\}$ is a linear function of
the variation~$\xi$, being the derivative of a smooth function.

Therefore we can collect these into a linear operator
$\Ag=\delta\dist$ from $\VF(L)$ to the space $C(\gStrut(L))$
of continuous functions on struts, defined by
$$ (\Ag\xi)(\{x,y\}) :=
 \delta_\xi \dist\{x,y\} =
 \frac{1}{|y-x|}\,\langle \xi_y - \xi_x, y - x \rangle. $$
Borrowing again from rigidity theory,
where the analogous~$A$ is called the \emph{rigidity matrix\/},
we will call~$\Ag$ the \emph{rigidity operator\/} for \gehring thickness.

The corollary above can be rephrased to conclude that a variation~$\xi$
decreases $\gThi(L)$ to first order if and only if~$\Ag\xi$ takes
at least one negative value on $\gStrut(L)$.

Note that, while the corollary says that \gehring thickness
has a directional derivative $\smash{\dplus_\xi} \gThi$ in each direction~$\xi$,
the operator $\smash{\dplus_\xi} \gThi$ is \emph{not\/} linear in~$\xi$.
For instance, when one component of a link is between two others,
it is easy to have both $\smash{\dplus_\xi} \gThi < 0$ and $\smash{\dplus_{-\xi}} \gThi < 0$.
We write the superscript~${}^+$ to emphasize that these are
only one-sided derivatives.  There is, however, a form of superlinearity:

\begin{corollary}\label{cor:dplus-super}
For any curve~$L$ and any $\xi,\eta\in\VF(L)$, we have
$$\dplus_{\xi+\eta}\gThi(L) \ge \dplus_\xi\gThi(L) + \dplus_\eta\gThi(L).$$
\end{corollary}
\begin{proof}
This follows immediately from the linearity of~$\Ag$ and the general
fact that $\min (f+g)\ge \min f + \min g$. We have
\begin{eqnarray*}
\min \Ag(\xi+\eta)\{x,y\} &=& \min (\Ag\xi + \Ag\eta)\{x,y\} \\
  &\ge& \min \Ag\xi\{x,y\} + \min \Ag\eta\{x,y\},
\end{eqnarray*}
where the minima are taken over all $\{x,y\}\in\gStrut(L)$.
\end{proof}

We will be interested in the adjoint~$\Ag^*$ of the rigidity operator,
so we first consider the dual function spaces.
By the Riesz representation theorem, we know that $C^*(\gStrut)$ is the
space of signed Radon measures on the space $\gStrut(L)$ of struts.
Similarly~$\VF^*(L)$ is the space of what we will call \emph{forces\/}
along~$L$, namely vector-valued Radon measures on~$L$.

The adjoint operator~$\Ag^*$ now associates to any measure~$\mu$ on struts
a force~$\Ag^*\mu$ along~$L$.  Geometrically, each pair $\{x,y\}$ acts
along the chord~$\overline{xy}$, outward at each of its endpoints.
In formulas,
\begin{eqnarray*}
\int_L \xi\,\d \Ag^*\mu &=& \int_{\gStrut} \!\!\Ag\xi\dmu \\
  &=&  \int_{x\in L}\! \int_{y\in L}\!
        \Big\langle\xi_y,\frac{y-x}{|y-x|}\Big\rangle\dmu(x,y),
\end{eqnarray*}
where we have lifted~$\mu$ to a symmetric measure $\mu(x,y)$ on ordered pairs.
Physically, we think of~$\mu$ as giving the strengths of compressive forces
within the struts, and~$\Ag^*$ as the operation that integrates
these strut forces to give a net force along the curve~$L$.

\section{First variation of length and finite total curvature}
\label{sec:ftc}

The objective functional we consider in this paper is simply
the length $\len(L)$ of a curve.  Since our curves might
not be smooth, we need to carefully examine the first variation
of length.

Let~$L$ be a rectifiable curve parametrized by arclength~$s$,
with unit tangent vector~$T$.
Suppose~$L_t$ is a variation of~$L$ under which the motion of
each point $x\in L$ is smooth in time with initial velocity~$\xi_x$,
and $\xi\in\VF(L)$ is a Lipschitz function of arclength.
Then the standard first-variation calculation shows that
$$\delta_\xi \len(L) := \ddtz{}{}{\len(L_t)} = \int_L \langle T,\xi'\rangle \ds,$$
where $\xi'=\d\xi/\d s$ is the arclength derivative, defined almost everywhere
along~$L$.

If~$L$ is smooth enough, we can integrate this by parts to get
$$\delta_\xi \len(L)
= -\int_L \langle T',\xi\rangle \ds - \sum_{x\in\bd L} \langle \pm T,\xi\rangle.$$
(In the boundary term, the sign is chosen to make~$\pm T$ point inward at~$x$.)
In fact, not much smoothness is required: 
as long as~$T$ is a function of bounded variation,
we can interpret~$T'$ as a measure, and the formula holds in a sense
we will now explore.

Following Milnor~\cite{milnor}, we recall that the total curvature of
a polygon is just the sum of its (exterior) turning angles, and we define
the \emph{total curvature\/} of any curve to be the supremal total curvature
over all inscribed polygons.  A rectifiable curve~$L$ has finite
total curvature if and only if the unit tangent vector $T=L'(s)$
is a function of bounded variation.  Sometimes the space of all such
curves is called $W^{1,\mathrm{BV}}$ or~$\mathrm{BV}^1$, but we will call it~$\FTC$.
(See Sullivan~\cite{Sul-FTC} for a survey of results on $\FTC$ curves.)

If $L\in\FTC$, it follows that at every point of~$L$ there are
well-defined left and right tangent vectors~$T_{\pm}$; these are equal
and opposite except at countably many points, the \emph{corners\/} of~$L$.
(See, for instance, Royden \cite[Sect.~5.2]{Royden}.)

If~$L$ is $\FTC$, its tangent~$T$ has a distributional derivative~$\K$
with respect to arclength: a force (an $\R^3$--valued Radon measure)
along~$L$ that we call the \emph{curvature force\/}.

The curvature force has an atom (a point mass or Dirac delta)
at each corner $x\in L$, with $\K\{x\}=T_+(x)+T_-(x)$.
On a $C^2$ arc of~$L$,  the curvature force is
$\K = \d T = \kappa N\ds$ and this is absolutely continuous
with respect to the arclength or Hausdorff measure $\ds=\H^1$.

When~$L$ has boundary, we choose to include in~$\K$ an atom
at each endpoint of~$L$, with mass~$1$ and
pointing in the inward tangent direction.  This means we
need no boundary terms in the formula
$\delta_\xi \len(L) = -\int_L \big<\xi,\d\K\rangle$.

We say that a vector field~$\xi$ along~$L$
is \emph{smooth\/} if~$\xi_s$ is a smooth function of arclength.
(The arclength parametrization of any rectifiable curve is
essentially unique, so this makes sense.)
The set of all smooth vector fields will be denoted $\VF^\infty(L)$.

The first variation $\delta \len(L)$ can be viewed as a linear
functional on smooth vector fields $\xi\in\VF^\infty(L)$.
As such a distribution, it has order zero, by definition, if 
$\delta_\xi \len(L)=\int_L \langle T,\xi'\rangle\ds$ is bounded by
$C\sup_L|\xi|$ for some constant~$C$.  This happens exactly when we
can perform the integration by parts.

We collect these results as:
\begin{lemma}\label{lem:ftc}
Given any rectifiable curve~$L$, the following conditions are equivalent:
\begin{enuma}
\item $L$ is $\FTC$.
\item The first variation $\delta \len(L)$ has distributional order zero.
\item There exists a curvature force measure~$\K$ along~$L$ such that
  $\delta_\xi \len(L) = - \int_L \langle \xi,\d\K\rangle$.
\end{enuma}
\end{lemma}

An $\FTC$ curve~$L$ is $C^1$ exactly when it has no corners,
that is, when~$\K$ has no atoms (except at the endpoints).
It is furthermore $C^{1,1}$ when~$T$ is Lipschitz,
or equivalently when~$\K$ is absolutely continuous (with respect
to arclength) and has bounded Radon--Nikodym derivative
$\d\K/\d\H^1 = \kappa N$.  In previous work on ropelength
(for example \cite{CKS2,MR2002m:74035}), the
thickness measure had an upper bound for curvature built in, meaning that any 
curve of positive thickness was automatically $C^{1,1}$. This
is not true for the \gehring thickness, so we do not expect
the same regularity results to hold here.

\section{Constrained criticality and the Kuhn--Tucker theorem}
\label{sec:constrainedcritical}

We will review constrained minimization problems in a finite setting,
before generalizing to the setting we will need for our ropelength problems.
Suppose we want to minimize a~$C^1$ function $f\co \R^n\to\R$
inside the \emph{admissible region\/} defined
by a finite collection of~$C^1$ inequality constraints~$g_i \ge 0$.
A constraint~$g_i$ is \emph{active\/} at~$p \in \R^n$ if $g_i(p) = 0$.

\begin{definition}
We say that~$p$ is a \emph{constrained critical point\/} for minimizing~$f$ if,
for any tangent vector~$v$ at~$p$ with $D_v f < 0$,
we have $D_v g_i < 0$ for some active~$g_i$.
That is, $p$ is critical if there is
no direction $v\in\R^n$ that reduces~$f$ to first order
while preserving all constraints to first order.
\end{definition}
Note that the criticality conditions for minimizing~$f$ and~$-f$
are quite different; in particular
local maxima for~$f$ are rarely critical points for minimizing~$f$,
while local minima for~$f$ usually---though not always---are.

\begin{figure}[ht!]\centering
\begin{overpic}[scale=.3]{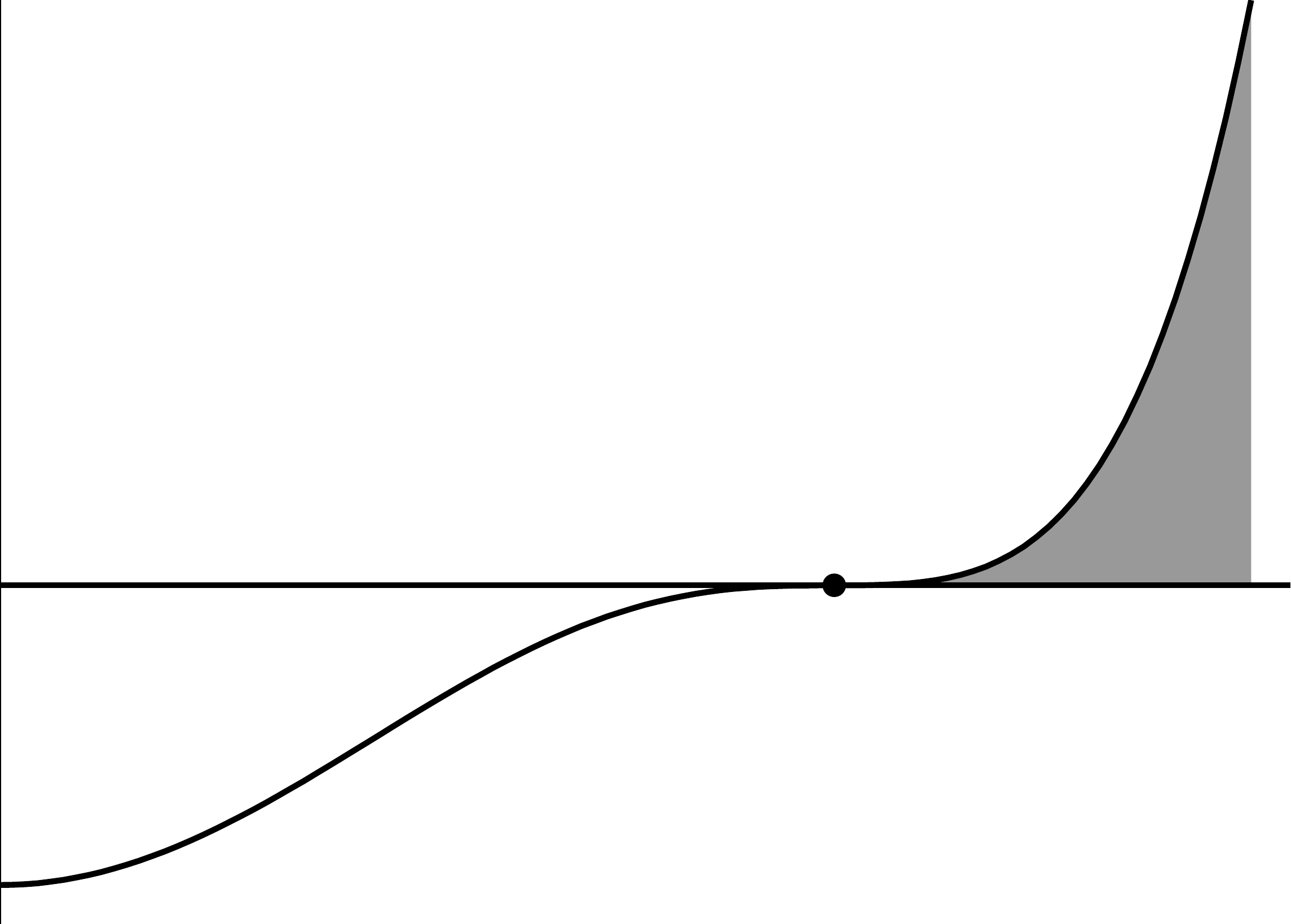}
\small
\put(11,31){$g_2=0$}
\put(76,62){$g_1=0$}
\put(62,21){$p$}
\put(96,22){$x$}
\put(1,66){$z$}
\end{overpic}
\caption[Example of unqualified constraints]
{In this illustration of \xmpl{min-not-crit},
the admissible region for the two constraints
$g_1=(x^2-1)^3-z\ge 0$ and $g_2=z\ge 0$ is shaded.
The Mangasarian--Fromovitz constraint qualification
fails at the cusp point $p=(1,0)$ because $\grad g_1$ and $\grad g_2$ are
equal and opposite there.}
\label{fig:unqualified}
\end{figure}
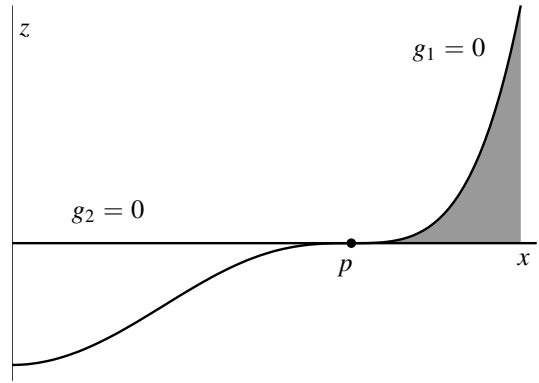

\begin{example}\label{ex:min-not-crit}
Suppose we minimize $f(x,z) := x$ on the halfplane $x\ge0$ in~$\R^2$, subject to
$$g_1 := (x^2-1)^3 - z\ge 0, \qquad g_2 := z\ge 0.$$
The admissible region has an outward-pointing cusp, shown
in \figr{unqualified}.
The tip of this cusp, at $p=(1,0)$, is the global minimum of~$f$
over the admissible region, but it is not critical:
the directional derivatives in the direction $v = (-1,0)$
are $D_v f = -1$ but $D_v g_i = 0$.
\end{example}

To deduce that a local minimum of~$f$ is critical according to our definition,
an additional regularity hypothesis will be required.
However, critical points can be exactly characterized
by a Lagrange multiplier theorem (compare \cite{MR13:855f}):
\begin{theorem}[Modified Kuhn--Tucker Theorem] \label{thm:ktfinite}
A point~$p$ is constrained-critical for minimizing~$f$ if and only if
the gradient~$\grad f$ is a positive linear combination
of the gradients~$\grad g_j$ of the constraints~$g_j$ active at~$p$.
\end{theorem}
The geometric intuition behind this theorem is easy to understand:
Only the active constraints matter, and being inequality constraints
they can only act positively.  If there were some component of $-\grad f$
not canceled by the $\grad g_j$, that would give an admissible direction to
move which decreases~$f$.

Unlike in the classical Kuhn--Tucker theorem, we do not need additional
regularity hypotheses on the point~$p$, which may surprise those
familiar with optimization theory. The explanation is that we are interested in
critical points, while the classical theorem deals with minima of~$f$.
And as we saw above, not every minimum of~$f$ is critical.
But just as in the classical theory, criticality will be guaranteed
if we add the hypothesis that
the Mangasarian--Fromovitz constraint qualification~\cite{MR34:7263}
holds for a local minimum.
\begin{definition} \label{defn:finitegregular}
A point~$p$ is \emph{constraint-qualified\/} (in the sense of
Mangasarian and Fromovitz) if there is a direction~$v$ such that
for all constraints~$g_j$ active at~$p$ we have $D_v g_j > 0$.
\end{definition}

We note that this condition fails
at the point $p=(1,0)$ in \xmpl{min-not-crit} above, which was
minimal but not critical.

\begin{proposition} \label{prop:finitelocalmincritical}
If~$p$ is a local minimum for~$f$ when constrained by $\{g_i\ge 0\}$,
and~$p$ is constraint-qualified,
then~$p$ is constrained-critical for minimizing~$f$.
\end{proposition}

We have omitted proofs of the theorem and proposition above
because they are standard and are also special cases of our
infinite-dimensional generalizations below.

\subsection{A generalized Kuhn--Tucker theorem}
Note that in \thm{ktfinite}, the functions~$f$ and~$g_i$
might as well be replaced by linear functions---their differentials at~$p$.
We view this as the linear-algebraic core of the Kuhn--Tucker theorem.

We will now derive an infinite-dimensional version,
where the linear functional~$f$ is defined on an arbitrary
vector space~$X$, and the finite family of constraints~$g_i$ is
replaced by a family~$A_y$, where~$y$ ranges over some compact space~$Y$.

While our theorem does not mention optimization directly,
it will be the engine that drives all of the
optimization theorems of this paper.

As usual, we let~$C(Y)$ be the Banach space
of continuous functions on~$Y$ with the sup norm $\|\cdot\|$,
and let $P\subset C(Y)$ be the closed positive orthant
consisting of nonnegative functions.
Then 
the dual space~$C^*(Y)$ consists of all signed Radon
measures on~$Y$, 
and $P^*\subset C^*(Y)$ is the cone of positive measures.

Note that any function $z\in C(Y)$ can be decomposed into positive
and negative parts: $z=z^+ - z^-$ with $ z^\pm \in P$.  Then
we have $\|z^-\|=\dist(z,P)$.

\begin{theorem} \label{thm:ourkt}
Let~$X$ be any vector space and~$Y$ be a compact topological space.
For any linear functional~$f$ on~$X$ and any linear map $A\co X \to C(Y)$,
the following are equivalent:
\begin{enuma}
\item\label{itm:strcrit} There exists $\eps>0$ such that
  $\|(A\xi)^-\| \ge \eps$ for all $\xi\in X$ with $f(\xi)=-1$.
\item\label{itm:balanced} There exists a positive Radon measure $\mu\in P^*$
  such that $f(\xi)=\mu(A\xi)$ for all $\xi\in X$.
\end{enuma}
\end{theorem}

This theorem is comparable to the generalized Kuhn--Tucker theorem of
Luenberger \cite[page~249]{luen}.  His theorem, restated to apply to
the linear Gateaux differentials ($f$ and~$A$) of the original objective
and constraint functions on~$X$, says:
\begin{theorem}[Luenberger]
Let~$X$ and~$Z$ be vector spaces, with a norm given on~$Z$,
and let $P\subset Z$ be a closed convex cone with nonempty interior.
Let $f\co X\to\R$ be a linear functional and $A\co X\to Z$ be a linear map.
Assume that whenever $A\xi\in P$ we have $f(\xi)\ge 0$,
and that~$A\xi$ lies in the interior of~$P$ for some $\xi\in X$.
Then there exists $\mu\in P^*$ such that $f(\xi)=\mu(A\xi)$
for all $\xi\in X$.
\end{theorem}

While our version applies only to the case $Z=C(Y)$,
our hypotheses \itm{strcrit} on~$f$ and~$A$ are somewhat weaker
than those imposed by Luenberger---they are
necessary as well as sufficient for \itm{balanced} the existence of~$\mu$.

To understand our overall strategy,
consider the linear map $(f,A) \co X \to \R \cross C(Y)$.
As we will see below, \itm{strcrit} implies that the image of $(f,A)$
avoids the interior of the orthant $\R^- \cross P$.

\begin{figure}[ht!]\centering
\begin{overpic}[scale=1.25]{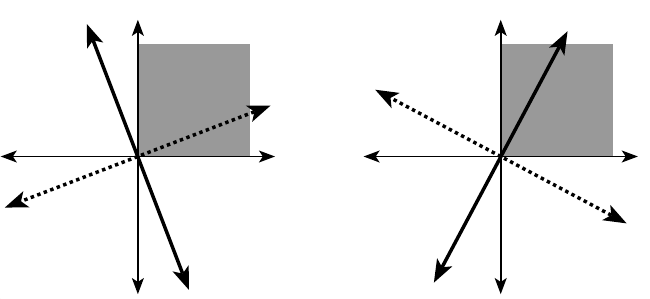}
\small
\put(39,40){$O$} \put(42,28){$I^\perp$} \put(10,41){$I$}
\put(95,40){$O$} \put(97,10){$I^\perp$} \put(85,42){$I$}
\end{overpic}
\caption[The Farkas alternative]
{One version of the Farkas alternative states that,
given any closed orthant~$O$ in a inner product space,
it must intersect at least one out of any
pair of orthogonal complements~$I$ and~$I^\perp$.}
\label{fig:farkas}
\end{figure}
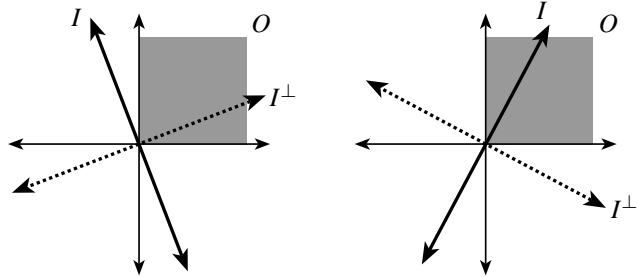

To gain some intuition, let us specialize to the case where $X= \R^m$.
  We can rephrase \itm{balanced} to say
that some vector in the kernel of the adjoint map $(f,A)^*$ is in $\R^- \cross P$.
When~$Y$ has finite cardinality~$n$, we put the standard
Euclidean inner product on $\R \cross C(Y)\isom \R^{n+1}$,
and identify this space with its dual.
Then the kernel of $(f,A)^*$ and the image of $(f,A)$
are orthogonal complements in~$\R^{n+1}$.
The standard Farkas alternative (see \figr{farkas}) says that,
given any closed orthant in~$\R^{n+1}$,
it must intersect one out of any pair of orthogonal complements.
Our argument in the general case, when~$Y$ might be infinite
and $X$ infinite dimensional, will be guided by this intuition. 

\begin{proof}[Proof of \thm{ourkt}]
One direction is easy:
suppose we have a positive Radon measure~$\mu$ so that for each $\xi \in X$,
\begin{equation*}
f(\xi) = \int_Y A\xi\dmu.
\end{equation*}
For any $\xi \in X$ with $f(\xi) = -1$, write $z:=A\xi\in C(Y)$.
We have $\int z\dmu=-1$, and since~$\mu$ is a positive measure,
we can replace the function~$z$ with its negative part
to conclude that $\int z^-\dmu \ge 1$.
Furthermore~$\mu$ has finite mass $\mass(\mu) := \int\d\mu < \infty$
by the Riesz theorem.
Therefore $$\dist(z,P)=\| z^- \| \ge 1/\mass(\mu)>0.$$
This completes the proof that \itm{balanced} implies \itm{strcrit}.

To prove the converse,
first give $\R \cross C(Y)$ the Euclidean combination of
the sup norms on~$\R$ and~$C(Y)$:
\begin{equation*}
\|(a,g)\| = \textstyle\sqrt{a^2 + \|g\|^2}.
\end{equation*}  
Now consider the orthant
$O := [-1,\infty)\times P$.
Our hypothesis \itm{strcrit} implies that there is positive distance
between~$O$ and the image $I:=(f,A)(X)$ of the linear map $(f,A)$.
Take sequences $\big(f(\xi_i),A\xi_i\big)$ in~$I$ and $(t_i,z_i)$ in~$O$,
whose pairwise distance approaches~$\dist(I,O)$.  That is, setting
$$v_i := \big(t_i-f(\xi_i), z_i-A\xi_i\big)$$
we have $\|v_i\|\approaches \dist(I,O)$.


We first claim that we can assume that $v_i\in \R^- \cross P$.
Certainly we can assume the equality $z_i=(A\xi_i)^+$, since this positive part of
the function~$A\xi_i$ realizes the distance 
$\| z_i - A\xi_i\| = \| (A\xi_i)^- \| = \dist(A\xi_i,P)$.
Then $z_i - A\xi_i = - (A\xi_i)^- \in P$.
Similarly, $\dist\big(f(\xi_i),\,[-1,\infty)\big)$
is $-1 -  f(\xi_i)$ if this is positive (and is zero otherwise),
so we may assume $t_i = \min(-1,f(\xi_i))$.
Thus $t_i \le f(\xi_i)$, so $t_i - f(\xi_i) \le 0$.
This proves the first claim.

We now have a geometric problem: from \figr{farkas} we see
there is a special case where both the image of $(f,A)$ and its orthogonal 
complement lie on the boundary of $\R^- \cross P$.  If this happens, then the
closure~$\bar I$ intersects either the
subspace $\R \cross \{0\}$ or $\{0\} \cross P$. 
The second case does not trouble us, but the first would cause us problems
later; we now show that our assumption \itm{strcrit} rules it out.
To do so, we think about the setup above
geometrically: if~$\bar I$ intersects $\R \cross \{0\}$, then 
we expect that $t_i - f(\xi_i) \rightarrow 0$. 

Thus our second claim is that we can assume
the $t_i - f(\xi_i)$ are uniformly negative.
If not, $\lim f(\xi_i) \le -1$, so without loss
of generality, we can rescale~$\xi_i$ down so that $f(\xi_i) = -1$.
That means $(-1,A\xi_i)\in I$.
By hypothesis~\itm{strcrit} we know
$$d_i := \dist\big((-1,A\xi_i),\,O\big)=\dist(A\xi_i,\,P) \ge \eps$$
for some fixed $\eps>0$. 
Since we are using the Euclidean combination of the norms on~$\R$ and~$C(Y)$,
the distance from any rescaling by~$k$ of~$(-1,A\xi_i)$ to~$O$ is given by
the Euclidean distance from~$(-k,k \| A\xi_i^-\|)$ to~$(-1,0)$.
And we can use plane geometry to see that rescaling by $1-\eps^2$
brings us closer to~$O$:
$$
\dist\big((1-\eps^2)(-1, A\xi_i),\, O\big)
= d_i \sqrt{1-2\eps^2+\eps^4(1+1/d_i^2)} 
\le d_i \sqrt{1-\eps^2+\eps^4}.
$$
We can always assume that $\eps < 1$, so the constant
$\sqrt{1-\eps^2+\eps^4}$ is less than~$1$. Therefore
\begin{eqnarray*}
\dist(I,O)
 &\le& \lim \dist\big((1-\eps^2)(-1, A\xi_i),\, O\big) \\
 &<& \lim \dist\big((-1,A\xi_i),\,O\big) \\
 &=& \dist(I,O).
\end{eqnarray*}
This contradiction proves the second claim.

We have proved that the~$v_i$ are in $\R^- \cross P$. 
Using the Hahn--Banach theorem, for each~$i$
we can find a linear functional $(c_i,\nu_i)\in\R \cross C^*(Y)$
that vanishes on~$I$, satisfies $(c_i,\nu_i)(v_i)=1$,
and has norm
$$\big\|(c_i,\nu_i)\big\|=1/\dist\big(v_i+I,\,(0,0)\big).$$
Because the $\R$--components of~$v_i$ are uniformly negative,
so are the~$c_i$.

Using Alaoglu's theorem,  the $(c_i,\nu_i)$ have a subsequence converging 
in the weak${}^*$ topology to a limit functional~$(c,\nu)$; we have $c < 0$ and
its norm is bounded above by $1/\lim\dist(v_i,-I)=1/\dist(I,O)$.

Setting $\mu:=\nu/|c|\in C^*(Y)$, we claim this
will be the Radon measure in statement~\itm{balanced}.
By construction, $(-1,\mu)$ vanishes on~$I$, meaning that
for $\xi \in X$, we have
$$ -f(\xi) + \int_Y A\xi\, \d\mu = 0.  $$
(Notice that we have used the additional geometric information
that~$I$ does not approach $\R \cross \{0\}$ in an essential way; if it did,
then~$c$ would vanish, and we could not rescale~$\nu$
by $1/|c|$ to obtain the equation above.)

It remains only to show that~$\mu$ is positive. In an inner product space,
this would be obvious: each~$\nu_i$ would be positive (since it was dual
to a positive $z_i - A\xi_i$), and~$\nu$ would be a limit of positive
measures. But our~$\nu_i$ were constructed implicitly by the Hahn--Banach
theorem, and so might include negative pieces. We now address this problem.

We can decompose each~$\nu_i$ into its positive and negative parts
$\nu_i=\nu_i^+ - \nu_i^-$,
with $\mass(\nu_i) = \mass(\nu_i^+) + \mass (\nu_i^-)$.
In order to show~$\nu$ is positive,
we will prove that $\lim \mass(\nu_i^+) = \lim\mass(\nu_i)$.
By construction, we know that
\begin{equation*}
1 = (c_i,\nu_i)(v_i) = c_i (t_i - f(\xi_i)) + \int_Y z_i - A\xi_i \,\d\nu_i.
\end{equation*}
Since $z_i - A\xi_i \in P$, we have
$$
\int_Y \! z_i-A\xi_i \,d\nu_i 
\le \int_Y \! z_i-A\xi_i \,d\nu_i^+
\le \big\|z_i - A\xi_i\big\| \mass(\nu_i^+).
$$
Using Cauchy--Schwarz, and the two equations above, we get
$$
1 \le \big\|v_i\big\| \,\sqrt{|c_i|^2 + (\mass(\nu_i^+))^2}.
$$
Now $\|v_i\|$ converges to $\dist(I,O)$,
so we find $\lim\|(c_i,\nu_i^+)\|\ge 1/\dist(I,O)$.
But the limit of $\|(c_i,\nu_i)\|$ (which cannot be smaller)
equals $1/\dist(I,O)$.
Therefore, $\lim\mass(\nu_i^+) = \lim\mass(\nu_i)$, completing the proof.
\end{proof}

To apply \thm{ourkt} to optimization problems, we will
let~$X$ be the space of variations~$\xi$ of our given configuration
and~$Y$ be the set of active constraints.
Then we let~$f(\xi)$ and~$A\xi(y)$ be the
directional derivatives of the objective function
and of the constraint $y\in Y$.

In this context, a configuration satisfying condition
\itm{strcrit} of \thm{ourkt} is called \emph{strongly critical\/},
and one satisfying \itm{balanced} is \emph{balanced\/}.
The theorem then says that a configuration is strongly critical
if and only if it is balanced.

Note that our strong criticality is indeed stronger than a simple
criticality condition, which would say that whenever $f(\xi)=-1$
we have $\dist(A\xi, P)>0$, or equivalently that no~$\xi$ has
$f(\xi)<0$ but $A\xi\in P$.
\begin{example}\label{ex:strcrit-not-crit}
With $X=\R^2$ and $Y=[0,1]$, we can set
$f(x_1,x_2)=x_1$ and $$A(x_1,x_2)(y)=2x_1 \tsty\sqrt{y-y^2} + x_2 y$$
to give an example that is critical,
but not strongly critical (and thus not balanced).
\end{example}

However, when~$Y$ is a finite set (with the discrete topology),
strong criticality is equivalent to criticality.  For suppose
whenever $f(\xi)=-1$ we have $\dist(A\xi, P)>0$, but there
is no uniform lower bound $\eps>0$ on this distance.
For each $y\in Y$, we know that $A\xi(y)$ is a linear functional on~$\xi$.
Since there are only finitely many~$y$, the graph of
$\min_{y\in Y} A\xi(y)$ describes a polyhedron in $X \cross \R$.
Since the supremum over $\xi\in X$ is finite (we know it is nonpositive),
it is achieved (at some~$\xi$ corresponding to a vertex of this polyhedron).
But for any~$\xi$, the value is negative, so this supremum must be negative.

This allows us to recover the finite-dimensional Kuhn--Tucker theorem:
let~$X$ be the tangent space to~$\R^n$ at~$p$,
let~$Y$ be the finite set of active constraints at~$p$,
and let~$f$ and~$A$ be the directional derivatives of
the objective function and the active constraints.
Because~$Y$ is finite, \itm{strcrit} is equivalent to
the definition of constrained criticality above,
and we obtain \thm{ktfinite}.

\section{The balance criterion for the Gehring problem} \label{sec:gehrbal}

We now have all the tools we need to develop a
balance criterion characterizing critical configurations
for the \gehring ropelength problem.
We start with definitions of criticality, guided by our
version of Kuhn--Tucker.

\begin{definition}
\label{defn:localmincritstrongcrit}
Suppose~$L$ is a rectifiable link with $\gThi(L)=\tau$,
and consider the Gehring problem of minimizing
length subject to the constraint $\gThi\ge \tau$. 
We say that~$L$ is:
\begin{itemize}
\item a \emph{local minimum\/} for \gehring ropelength
if for all~$L'$ sufficiently $C^0$--close to~$L$ with $\gThi(L')\ge \tau$
we have $\len(L')\ge\len(L)$.
\item \emph{critical\/}
if for all $\xi \in \VF^\infty(L)$
with $\delta_\xi \len(L) < 0$ we have $\dplus_\xi \gThi(L) < 0$.
\item \emph{strongly critical\/}
if there exists some $\eps > 0$ such that
for all $\xi \in \VF^\infty(L)$ with $\delta_\xi \len(L) = -1$,
we have $\dplus_\xi \gThi(L) \le -\eps$.
\end{itemize}
\end{definition}

With these definitions,
we can now apply our Kuhn--Tucker theorem to the Gehring problem.

\begin{theorem}[Balance Criterion] \label{thm:gbalance}
A link~$L$ is strongly critical for length when
constrained by \gehring thickness if and only if there exists a positive
Radon measure~$\mu$ on $\gStrut(L)$ such that,
for every smooth vector field~$\xi$ along~$L$, we have
$$\delta_\xi\len(L) = \int_{\gStrut(L)}\Ag\xi\dmu,$$
where $\Ag=\delta\dist$ is the rigidity operator.
\end{theorem}
\begin{proof}
We will apply \thm{ourkt} with $X:=\VF^\infty(L)$ and $Y:=\gStrut(L)$,
letting $f:=\delta\len(L)$ be the derivative of length and $A:=\Ag$
be the rigidity operator.  We have
$$\|(\Ag\xi)^-\| = -\min_{\gStrut} \delta_\xi \dist\{x,y\}$$
(when this is nonnegative).
By \cor{dgthi}, the right-hand side is $-\dplus_\xi\gThi(L)$,
so that condition \itm{strcrit} from \thm{ourkt}
is exactly strong criticality.
\end{proof}

\subsection{Smoothness of critical curves}
It is unclear, a priori, how much regularity one should expect
for ropelength-critical curves in the Gehring problem.  But we can use
the balance criterion to deduce immediately that they must
have finite total curvature.

\begin{corollary}\label{cor:gehr-ftc}
If a link~$L$ is strongly critical for the Gehring problem,
then~$L$ is~$\FTC$.
\end{corollary}
\begin{proof}
The theorem tells us that~$L$ can be balanced:
$$\delta_\xi\len(L) = \int_{\gStrut(L)}\!\!\Ag\xi\dmu.$$
But the right-hand side is a distribution of order zero on~$\xi$, since
$$\int_{\gStrut(L)}\!\!\Ag\xi\dmu \le \mass(\mu)\,\sup_L|\xi|.$$
Therefore, by \lem{ftc}, $L\in\FTC$.
\end{proof}

We can now rewrite the conclusion of our balance criterion in
terms of the curvature force~$\K$ on~$L$ and the adjoint~$\Ag^*$
of the rigidity operator.

\begin{corollary}\label{cor:gehr-bal}
A link~$L$ is strongly critical for \gehring ropelength if
and only if it has finite total curvature and there exists
a positive Radon measure~$\mu$ on $\gStrut(L)$ such that
$$\Ag^*(\mu)=-\K$$
as forces along~$L$.
\end{corollary}
\begin{proof}
The theorem guarantees that for all smooth~$\xi$, we have
$$\delta_\xi\len(L) = \int_{\gStrut(L)}\Ag\xi\dmu.$$
By the corollary, $L$ is $\FTC$, so the left-hand side
can be rewritten as $-\int_L \langle\xi,\d\K\rangle$.  Approximating
any continuous vector field uniformly by smooth ones, we find that
$$-\int_L\langle \xi,\d\K\rangle= \int_{\gStrut(L)}\Ag\xi\dmu$$
for all $\xi\in\VF(L)$, or in other words, $-\K=\Ag^*(\mu)$.
\end{proof}

We get an immediate and useful geometric corollary to this balance criterion.

\begin{corollary}\label{cor:strut-cvxhl}
Suppose~$L$ is critical for \gehring ropelength,
and $E\subset L$ is a subset with nonzero net (vector) curvature
$0\ne\K(E)\in\R^3$.
Then there must be at least one strut $\{e,x\}$ with $e\in E$ and
$x\notin E$, and~$\K(E)$ is in the convex cone generated by
the directions $x-e$ of all such struts.
\end{corollary}
\begin{proof}
First note that struts from~$E$ to~$E$ contribute no net force.
By the balance criterion, we have $\K(E)=-\Ag^*\mu(E)$, and the
latter is a (positive) weighted sum of vectors $x-e$.
\end{proof}
We note that this corollary is the analogue for \gehring ropelength
of von der Mosel and Schuricht's
``Characterization of Ideal Knots''~\cite[Theorem~1]{vdms}.

We next find that critical links are~$C^1$ as well as $\FTC$:
\begin{proposition} \label{prop:gcrit-c1}
If~$L$ is strongly critical for \gehring ropelength, then~$L$ is~$C^1$.
\end{proposition}

\begin{figure}[ht!]\centering
\begin{overpic}[scale=1.5]{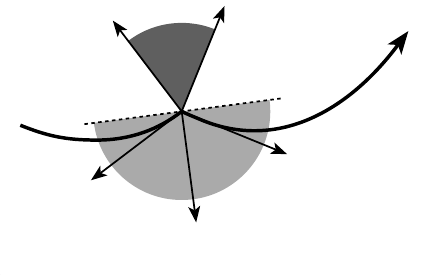}
\small
\put(16,18){$T_{-}$}
\put(39,7){$\K(\{x\})$}
\put(69,25){$T_{+}$}
\put(91,44){$L$}
\put(45,40){$x$}
\end{overpic}
\caption[\gehring critical curves have no corners]
{This curve~$L$ has a corner at~$x$ with left
and right tangent vectors~$T_{-}$ and~$T_{+}$, whose sum is the curvature
force~$\K(\{x\})$ there. If~$L$ is to balance, there must be a strut $\{x,y\}$
with~$y$ in the open hemisphere (shown in light gray)
of vectors with positive inner product with~$\K(\{x\})$.
But for any~$y$ outside the \emph{normal cone\/} (shown in dark gray),
there are points near~$x$ on~$L$ that are closer to~$y$ than~$x$ is.
Thus our $\{x,y\}$ cannot be a local minimum of the self-distance function.
This contradiction proves that a critical curve cannot have a corner.}
\label{fig:gehr-corn}
\end{figure}
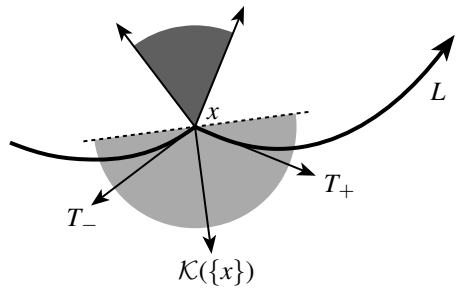

\begin{proof}
We already know that~$L$ has finite total curvature;
it is~$C^1$ precisely when it has no corners, that is,
when the curvature force~$\K$ has no atoms.
If~$T_\pm$ are the right and left tangent vectors to~$L$ at~$x$,
then $\K(\{x\}) = T_+ + T_-$.  When $\K(\{x\}) \neq 0$,
\cor{strut-cvxhl} says there is
at least one strut $\{x,y\}$ with $\langle y-x, \K(\{x\}) \rangle >0$.
That is,
$$ \langle y-x,  T_+\rangle + \langle y-x, T_-\rangle > 0,$$
so we must have $\langle y-x, T_+\rangle >0$ or $ \langle y-x, T_-\rangle >0$.
(See \figr{gehr-corn}.)
In either case it follows that there exist points on~$L$ near~$x$
that are closer to~$y$ than~$x$ is, which contradicts the hypothesis
that $\{x,y\}$ was a strut. This completes the proof.
\end{proof}

The example of the tight clasp in \secn{gehr-clasp} shows that
critical links need not be $C^{1,1}$---their
curvature need not be bounded---but so far this is the worst behavior we
can display. We conjecture that the curvature measure is always
absolutely continuous with respect to arclength.

\subsection{Constraint qualification in the sense of Mangasarian--Fromovitz}
\cor{gehr-bal} will be the basic model for balance criteria
for generalized links, and for links constrained by
other thickness functionals~\cite{CFKSW2}.
In some cases, including the \gehring ropelength for closed
links we are treating now, we can
improve on this form of the criterion by replacing strong criticality with
criticality. This is our next goal.

In \secn{constrainedcritical}, we defined a
regular or constraint-qualified point
for a finite set of constraints $g_1, \dots, g_n$:
such a point has some variation direction~$v$ such that $D_v g_i > 0$
for all the active~$g_i$.  By \cor{dgthi},
the corresponding idea for a link~$L$ in the Gehring problem
is the existence of a vector field~$\xi$ for which $\dplus_\xi \gThi(L) > 0$.
But this is automatic: dilating~$L$ increases $\gThi$ to first order.

This regularity for our problem allows us to prove that local
minima are critical and that critical points are strongly critical.
\begin{proposition}\label{prop:gmincrit}
A link~$L$ is critical for the \gehring ropelength problem if
and only if it is strongly critical.
If~$L$ is a local minimum, then~$L$ is critical.
\end{proposition}
\begin{proof}
Suppose~$L$ is a local minimum but not critical.
Then for some $\xi\in \VF^\infty(L)$ we have
$\delta_\xi\len(L)<0$ but $\smash{\dplus_\xi}\gThi(L)\ge 0$.
Then for small enough $t>0$,
the link $L+t\xi$ has less length than~$L$.
This contradicts minimality unless $\delta_\xi\gThi(L)=0$ and thickness
has decreased (but not to first order).  But in this case, we can
instead use the rescaled deformation $(\gThi(L+t\xi))^{-1}(L+t\xi)$, for which $\gThi \equiv 1$.
For small $t>0$ these again have less length than~$L$, contradicting
minimality.

Strong criticality always implies criticality.
Conversely, suppose a closed link~$L$ is critical but not strongly critical.
Then there exists a sequence
$\xi_i \in \VF(L)$ with $\smash{\delta_{\xi_i}} \len(L) = -1$ and
$\smash{\dplus_{\xi_i}} \gThi(L) \rightarrow 0$.
Let~$\eta$ be the vector field along~$L$ induced by dilation,
scaled so that $\delta_\eta \len(L) < 1$.
Then we observe that
$\delta_{\eta + \xi_i} \len(L) < 0$
for all~$i$.
The superlinearity of \cor{dplus-super} shows that
$$\lim\dplus_{\eta + \xi_i} \gThi \ge \dplus_{\eta} \gThi > 0.$$
But then for some~$i$, we must have 
$$\delta_{\eta+\xi_i} \len(L) < 0 \quad \text{and} \quad \dplus_{\eta+\xi_i} \gThi(L) > 0,$$
contradicting the criticality of~$L$.
\end{proof}

Thus for closed links, a minimizer (or more generally any critical point) for
the \gehring ropelength problem is strongly critical, and hence by
\cor{gehr-bal} its curvature force is balanced by some strut force~$\Ag^*\mu$.
However, in our generalized ropelength problems, with endpoint constraints
and obstacles, constraint qualification will not always hold.
Then we will have to be careful 
about the distinction between criticality and strong criticality.

\subsection{Existence of minimizers}
We now show that each link-homotopy class
contains a globally length-minimizing curve with $\gThi \ge 1$.
\begin{proposition}
\label{prop:gminexist}
In a given link-homotopy class~$[[L]]$,
among all curves with \gehring thickness at least~$1$,
there is some~$L_0$ of minimum length.
\end{proposition}

\begin{proof} 
We may rescale the initial~$L$ so that $\gThi(L) \ge 2$.
Thus if the $C^0$ distance between~$L$ and a link~$L'$
is less than $\nicefrac{1}{2}$
then the straight-line homotopy between them is a link homotopy,
and $\gThi(L')\ge1$.
Taking $L'$ to be a standard smoothing of~$L$
(for instance, its convolution with a smooth bump function),
it follows that $[[L]]$ contains a $C^\infty$ link.

In particular, the set of rectifiable links in $[[L]]$ with
link-thickness at least~$1$ is nonempty. Let $L_1,L_2,\dots$ be a sequence
of such links with lengths tending to the infimal length~$\ell$ in
this class. By the Arzela--Ascoli theorem, taking a subsequence we may
assume that the~$L_i$ converge in~$C^0$ to a limit~$L_0$.  Since
$\gThi$ is continuous with respect to the $C^0$ topology, and length
is lower semicontinuous, it follows that $\gThi(L_0)\ge 1$ and
$\Len(L_0)\le \ell$.  By the remarks of the last paragraph, $L_0$ is link
homotopic to~$L_i$ for large~$i$, and therefore $L_0 \in
[[L]]$. Thus $\Len(L_0) = \ell$ and $L_0$ is the required minimizer.
\end{proof}

Since $C^\infty$ links are tame, the argument above also
shows the following:

\begin{proposition} \label{prop:tame}
There are no wild link-homotopy types with finitely many components.
\end{proposition}

(This was originally observed by Milnor~\cite{MR17:70e}.)
Thus in the work to come, we need only to consider tame links.

\section{Examples of critical links} \label{sec:gehr-examples}

\subsection{The known length-minimizing links}
In~\cite{CKS2},
we showed that if one component of a link is linked to~$k$ others
then its length is at least a certain constant~$P_k$.
Although our theorem was written for the original ropelength problem,
the proof is valid for the Gehring problem as well.
Whenever a link can be realized with each component having length~$P_k$,
that configuration is thus a length-minimizer
not only when constrained by thickness but also when
constrained by \gehring thickness.
(These are still the only examples known to be ropelength-minimizers.)

To any link~$L$ we can associate a graph: the vertices are
the components of~$L$, and the edges record which pairs are
nontrivially linked.  For any tree~$T$ with~$n$ edges, there is
a unique link $H(T)$ that is a connect sum of~$n$ Hopf links
and whose associated graph is~$T$.

\begin{figure}[ht!]
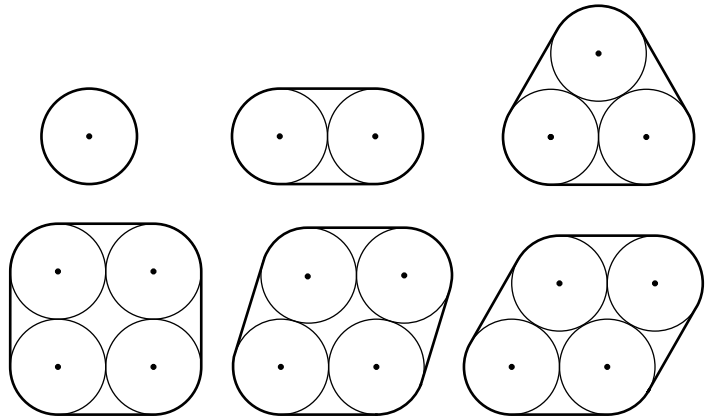

\centerline{\hspace{3mm} \incgr{1.0}{\figdir/peri-1} \hspace{10mm} \incgr{1.0}{\figdir/peri-2} \hspace{8mm} \incgr{1.0}{\figdir/peri-3} } 
\vspace{12pt}
\centerline{\incgr{1.0}{\figdir/peri-4}}
\caption[Perimeter-minimizing enclosures of disks in the plane]
{Here we see perimeter-minimizing enclosures of $n=1$,
$2$,~$3$ and~$4$ unit-diameter disks in the plane.  The components
in the known minimizing links are outer parallels to such
curves at distance~$\nicefrac{1}{2}$.  When $n=4$, the minimizer
does not have a unique shape; instead there is a one-parameter
family of minimizers.  In the last shape on the lower left,
there is one additional isolated strut, but it carries no
force in the balancing measure.}
\label{fig:peritable}
\end{figure}

For many trees~$T$ with vertices of sufficiently low degree,
we can realize $[[H(T)]]$ explicitly
with each component having exactly its minimum possible length~$P_k$.
Even some slightly more complicated links, like the example in \figr{known},
whose graph is not a tree, can be realized in this way.
\figs{known}{1.4}{\figdir/newcomplex}
{A minimizing link for the \gehring ropelength problem.}
{This link of six components is a global minimizer for
the \gehring ropelength problem.  Each component is a convex plane
curve that minimizes its length given the number of other components
it links.}
The distance between any two linked components is exactly~$1$.
Each component in one of these minimizers is a convex plane curve
built from circular arcs of radius~$1$ and straight segments.
It is an outer parallel (at distance~$\nicefrac{1}{2}$) to a shortest
curve surrounding~$n$ disjoint unit-diameter disks in the plane.
(See \figr{peritable}.)

Consider the $n$--star~$T_n$, the tree with a central vertex incident
to all~$n$ edges.  For $n\le5$, the construction above produces a
\gehring ropelength-critical configuration of~$H(T_n)$ that is known to be minimizing.
We will examine the case $n=2$ in detail, in light of our balance criterion,
and then indicate how to produce \gehring ropelength-critical configurations for all~$n$.

\begin{example}\label{ex:hopf-chain}
The link $H(T_2)$ is the simple chain of three components,
shown in \figr{chainstruts}.
In the ropelength minimizer, the two end components are circles~$C_1$ and~$C_2$,
while the middle component is a stadium curve~$S$.  The centers
of the circular arcs in~$S$ are points $c_i\in C_i$, while the
center of each~$C_i$ is a point $s_i\in S$.
The struts are exactly where different components
are at distance~$1$.  There is a strut from each point along
each circular arc to the center of that arc (from~$C_i$ to~$s_i$
and from~$S$ to~$c_i$).  There is also one further strut $\{c_1,c_2\}$.

\begin{figure}[ht!]\centering
\begin{overpic}[scale=1.25]{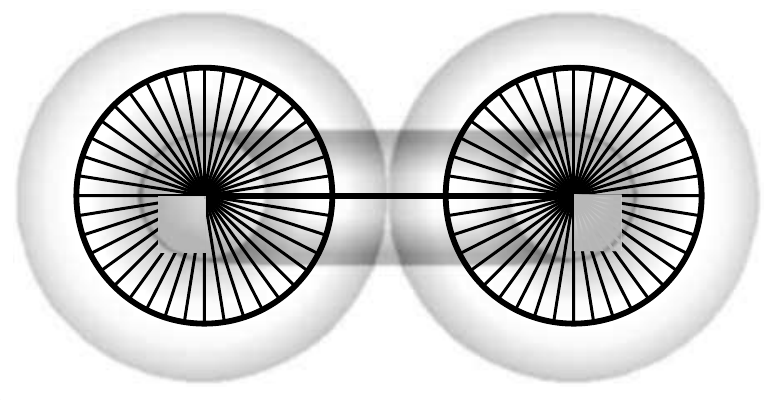}
\small
\put(22,22){$s_1$}
\put(11,10){$C_1$}
\put(85,10){$C_2$}
\put(75.5,22){$s_2$}
\put(48,22){$S$}
\end{overpic} \hfil
\begin{overpic}[scale=1.25]{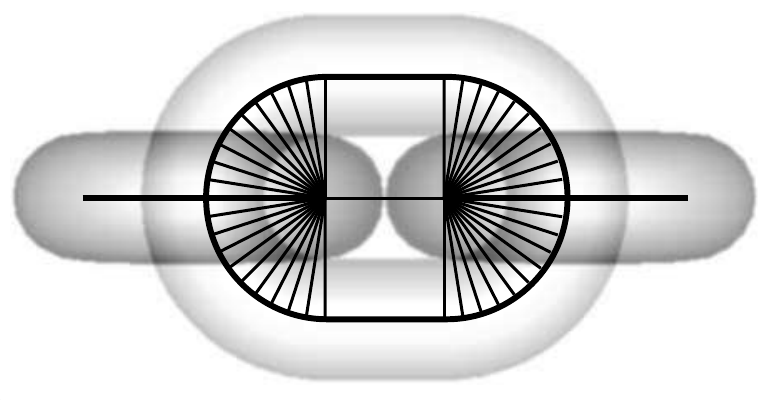}
\small
\put(43,23){$c_1$}
\put(53,23){$c_2$}
\put(48,44){$S$}
\put(11,22){$C_1$}
\put(85,22){$C_2$}
\end{overpic}
\caption[Struts to balance the simple chain]
{This simple chain is known to be a minimizer for the \gehring ropelength
problem, so by the balance criterion, its curvature force must be balanced
by some measure on the struts.
At the top, we see how the curvature forces along the circular
components~$C_i$ are balanced by the struts coming into the centers~$s_i$.
They produce no net force on either center~$s_i$.  At the bottom, we see
how the curvature forces along the semicircles of~$S$ are balanced by
struts to their centers~$c_i$.  The resulting net inward force on the~$c_i$
is balanced by an atomic measure on the one remaining strut $\{c_1,c_2\}$.}
\label{fig:chainstruts}
\end{figure}

Since we know that this configuration is
length-minimizing when constrained by \gehring thickness,
these struts, by \cor{gehr-bal} and \prop{gmincrit},
must support a balancing measure~$\mu$.
Conversely, exhibiting such a measure will re-prove that this configuration
is critical for the \gehring ropelength problem, though
to re-prove it is a local minimum would require some second-order theory.
We now provide such a measure, which will be a useful comparison of the results
of this paper against the results of~\cite{CKS2}.

Except for~$c_i$, each point~$x$ along the component~$C_i$ is part of
a unique strut $\{x,s_i\}$.  The measure assigned to struts in this ``wheel''
must exactly balance the curvature force $\d\K=N\ds(x)$ along~$C_i$.
Because the wheel forms a complete circle,
at the center points~$s_i$, the incoming forces from these struts
cancel one another, leaving no net force.

The situation on the stadium curve is slightly more complex. 
The struts from the semicircles of~$S$ to the points~$c_i$ again
balance $\d\K=N\ds(x)$, now for~$x$ along the semicircles.
Unlike the previous situation, however, these measures have
a resultant inward force of magnitude~$2$ at~$c_i$,
directed parallel to the straight segments in the stadium curve.
To balance these forces, the measure~$\mu$ must have an atom
of magnitude~$2$ at the one remaining strut $\{c_1,c_2\}$.

The measure~$\mu$ we have described does balance the curvature force
everywhere along the link, and thus demonstrates that the link is critical
for \gehring ropelength.

It is worth emphasizing the fact that the inner strut $\{c_1,c_2\}$
bears an atom of~$\mu$.  This stresses the point that in our Kuhn--Tucker
theorem and the resulting balance criterion we are required to 
view the Lagrange multiplier~$\mu$ as a Radon measure in the dual
space $C^*(\gStrut(L))$, rather than as a density function on struts.
\end{example}

Although ropelength-minimizing, \xmpl{hopf-chain}
is not rigid, in the sense that the components~$C_i$ can be pivoted
around the points~$c_i$ to be centered at any points~$s_i$ on
the semicircles of~$S$.

A stronger form of nonuniqueness is exhibited by
the minimizing configurations~\cite{CKS2} of the five-component
link $H(T_4)$, with one component linked to all four others.
Here the central component does not even have a uniquely
determined shape.  Instead there is a one-parameter family of
minimizing shapes, corresponding to the deformation seen
in \figr{peritable} for $n=4$.
Again, each of the minimizers can be balanced.
(As we have proven, the existence of the balancing measure~$\mu$
is equivalent to strong criticality for the ropelength problem,
but it does not imply that the critical point is isolated.)

For $n>5$, we expect that similar configurations of $H(T_n)$,
like the one shown in \figr{wrapped} for $n=6$,
are again minimizing.  Our balance criterion lets us show they
are at least critical:
\begin{figure}[ht!]\centering
\begin{overpic}[scale=.45]{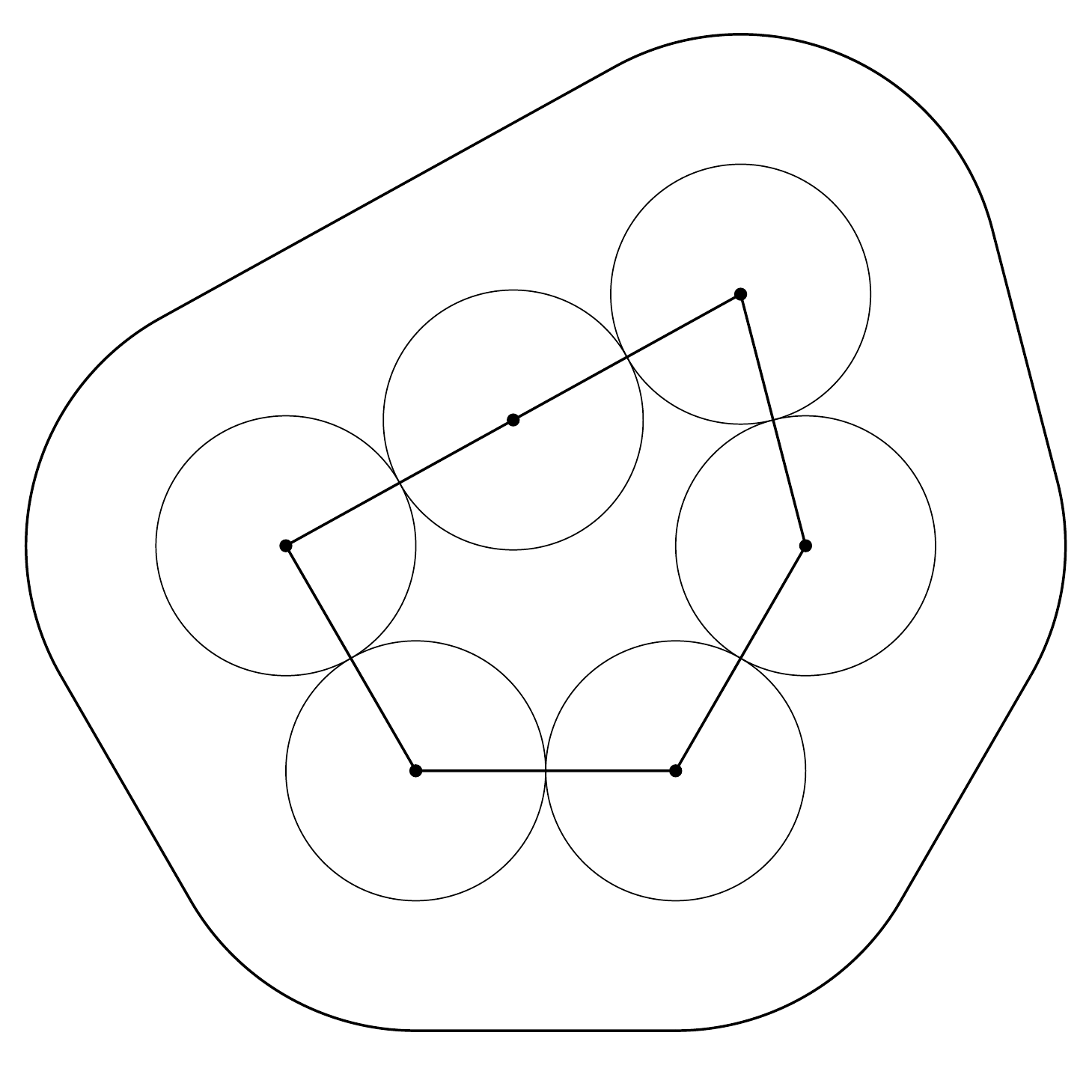}
\small
\put(7,18){$L_0$}
\put(70,66){$P$}
\end{overpic}
\caption[Wrapped link]
{This configuration of $H(T_6)$ is critical
for the \gehring ropelength problem.  It is also presumably
the minimizer, even though it does not minimize the
length of the long component~$L_0$ alone.  That component is the outer parallel
at distance~$1$ to a convex planar polygon~$P$.  Each other component~$L_i$
is a unit circle passing through a vertex of~$P$ and lying in a perpendicular
plane.  We have drawn the unit-diameter
disks around these vertices, where the thick tubes around the~$L_i$ intersect
the plane inside~$L_0$.}
\label{fig:wrapped}
\end{figure}

\begin{proposition}
\label{prop:wrapped}
Suppose~$P$ is a convex planar $n$--gon with unit-length sides 
and turning angles in $[0,\nicefrac{2\pi}{3}]$.  Let~$L_0$ be the outer parallel
at distance~$1$ from~$P$, and let $L_1,\ldots,L_n$ be unit circles,
perpendicular to the plane of~$P$, passing through the vertices of~$P$,
and centered at points on~$L_0$.
Then the link $L=L_0\cup L_1 \cup \cdots\cup L_n$
is a configuration of $H(T_n)$ with \gehring thickness~$1$
that is critical for \gehring ropelength. 
\end{proposition}

\begin{proof}
As in the simple chain, each circle~$L_i$ focuses a wheel of struts to its
center point on~$L_0$, and a measure assigning force~$\ds$ to these 
struts balances the curvature force on each circle while exerting no net force
on~$L_0$.

Let~$c_i$ be the vertex of~$P$ on~$L_i$, and let~$2\alpha_i$ be the
turning angle of~$P$ there.
The condition $\alpha_i\le\nicefrac{\pi}{3}$ exactly
suffices to know that no two vertices (and thus no two~$L_i$)
are at distance less than~$1$ from each other, confirming that $\gThi(L)=1$.
The curve~$L_0$ includes an arc of the unit circle around~$c_i$;
from this arc of length~$2\alpha_i$ a fan of struts converge to~$c_i$.
To balance the curvature force on~$L_0$, these struts again have measure
equal to~$\ds$, giving a net inward force of $2\sin\alpha_i$ on~$c_i$.
The remaining, isolated struts of~$L$
connect successive~$c_i$ along the edges of~$P$.
Unit atoms of compressive force on these isolated struts produce
exactly the outward forces $2\sin\alpha_i$ at~$c_i$ needed to
balance the inward forces from~$L_0$.

By \cor{gehr-bal}, the existence of this balancing measure on the
struts proves that~$L$ is critical.
\end{proof}

For $n\le 5$, we know these configurations for $H(T_n)$
are ropelength minimizers.
For $n>5$, the component~$L_0$, having length $n+2\pi$, is longer than
it needs to be: at the expense of lengthening some other components,
it could be shortened to length~$P_n$, which, asymptotically,
is much smaller, being $O(\sqrt n)$.
However, calculations we have done suggest that
the tradeoff is not worthwhile and so the critical configuration
described above is probably the global minimum for ropelength.

The examples given in \prop{wrapped}---critical configurations
and presumed minimizers for $H(T_n)$---are quite interesting.
The shape of~$L_0$ is free to move in an $(n-3)$--parameter family;
each other component is free to pivot (about its vertex of~$P$ and
along one of the arcs of~$L_0$),
giving an additional~$n$ parameters for the shape of the whole link~$L$.
We also note that these examples are tight links that are not
packed tightly:  Consider the thick (unit diameter) tube around one
of these configurations.
As~$n$ increases, it occupies an ever smaller fraction of the volume of
its convex hull.  This should be compared with
experiments of Millett and Rawdon~\cite{mrjcp} on this volume fraction.

Although we have stated \prop{wrapped} above only for stars~$T_n$, the
same balancing works for the links $H(T)$ based on other trees~$T$.
Each component linked to~$n$ others should have the shape of~$L_0$ above.
We note, however, that critical links built in this way are not always
minimizers.
\begin{example}\label{ex:loose-tree}
Consider the tree $T_{n,m}$ with $n+m$ vertices, including one of
valence~$n$ connected to another of valence~$m$.  The link $H(T_{n,m})$
then has two long components, $L_0$ and~$L_1$, linked to each
other and to $n-1$ and $m-1$ short components, respectively.
For large enough $n$ and~$m$, we can construct two thick versions
of $H(T_{n,m})$, called the~$A$ and~$B$ configurations, as follows.

For the $A$ configuration, we follow \prop{wrapped} to build~$L_0$
and~$L_1$ as outer parallels to convex planar $n$- and $m$--gons $P_0$
and~$P_1$ in perpendicular planes. If $n = 2k$, then we let~$P_0$ take
the shape of a $(k-1)\times 1$ rectangle capped with equilateral
triangles on both short sides; if $n = 2k + 1$, we
omit one triangle.  (The precise shape is unimportant, but we need
at least one sharp angle on each polygon.)  Further, we choose
the tip of such a triangle as the vertex of~$P_0$ corresponding to~$L_1$
(and as the vertex of~$P_1$ corresponding to~$L_0$).
Then it is straightforward
to check that the other components stay sufficiently far from each
other for this $A$ configuration to indeed have $\gThi=1$;
its total length $(2\pi+1)(n+m)$.

However, when~$n$ and~$m$ are large enough,
we can save two units of length as follows.
Construct a tight configuration of $H(T_{n-1})$
as in \prop{wrapped}, using a regular polygon.  (Again, the precise shape
is not important, but here we need a large hole in the middle of the polygon.)
This configuration, and indeed the unit-diameter thick tubes around
its components, is contained in a round solid torus $U_n$
of minor radius~$\nicefrac32$
and major radius $1+1/(2\sin\nicefrac{\pi}{n-1})$.
Then construct the analogous configuration of~$H(T_{m-1})$
contained in a solid torus $U_m$.  Finally, place these
two pieces in space so that~$U_n$ and~$U_m$
form a (loose) Hopf link.  (This is possible as long as the major
radii are at least~$3$, corresponding to $n,m\ge14$.)
The resulting link is the $B$ configuration of~$H(T_{n,m})$.
Because the large Hopf link is loose, there are no struts from~$H(T_{n-1})$
to~$H(T_{m-1})$.  Since each of these pieces is balanced, so is
the $B$ configuration of~$H(T_{n,m})$.

If, as we believe, these $B$ configurations are the ropelength minimizers,
then they are the first ones known in which certain pairs
of linked components are not in contact.  (We note that the
same must be true for the $n$--component Hopf links for large~$n$,
since their minimum ropelength~\cite{CKS1} is $O(n^{\nicefrac32})$.
There, however, no explicit candidate minimizer is known.  And
our $B$ configuration here has the additional property that certain
linked pairs are not even connected by chains of touching components.)
\end{example}

In all of the examples $H(T)$ discussed above, each component
is a convex plane curve built from straight segments and arcs of unit circles.
The proven minimizers are minimizers in their
isotopy class for the original ropelength problem~\cite{CKS2}, as well as
minimizers in their link-homotopy class for \gehring ropelength.
In fact, in~\cite{CFKSW2} we will consider a family of thicknesses with
varying stiffness.  Each of these thicknesses is characterized by a 
stiffness~$\lambda$, meaning a lower bound on the diameter of
curvature for a unit-thickness curve.
Our~$H(T)$ are ropelength-minimizers for the whole family, as long as the
stiffness~$\lambda$ does not exceed~$2$,
when circular arcs of larger diameter would be needed.  We will also develop
an analog of our balance criterion for these other ropelength problems,
and will see that all the~$H(T)$ discussed above (including those that
are not minimizers) are critical for all formulations of ropelength
where $\lambda\le2$.

\subsection{Local minima for ropelength}
We do not attempt in the paper to discuss second-order behavior of
ropelength near a critical point---in particular we have no way yet
to distinguish between local minima and saddle points for this problem.
Of course, the known minimizers must be local minima, and it is also
easy to give critical configurations which are not local minima,
as in \xmpl{unstable} below.

Many researchers have used numerical simulations of the ordinary
ropelength problem to look for nontrivial local minima for knots,
in particular for the unknot.  Such configurations have been termed
\emph{Gordian unknots\/} since they can be untangled topologically
but not physically.  Pieransky et al~\cite{Pier-gord}
have numerically simulated a reasonable candidate for a Gordian
unknot, but we are very far from being able to prove its existence.

In \xmpl{loose-tree} we gave two distinct critical configurations for
$H(T_{n,m})$, and we expect that this will lead to the provable
existence of two distinct local minima.  In particular, our investigations
lead us to predict that one cannot move from the $A$~configuration
of length $(2\pi+1)(n+m)$ to the suspected global minimum~$B$ without
first increasing ropelength.  This shows there must be a second
local minimum; we expect, however, that this is not~$A$ but instead
a third configuration of intermediate length.

This connects back to Alexander Nabutovsky's original work on ropelength
in higher dimensions and codimensions~\cite{Nabut}.  He showed
using recursive function theory that, in those higher dimensions,
a ropelength constraint often introduces new components into the
moduli space of unknotted hyperspheres; in particular there are
infinitely many local minima for ropelength.  While for two-spheres
in~$\R^3$ or for circles in~$\R^2$ there are presumably no such minima, we
do expect there must be infinitely many Gordian unknots in~$\R^3$.
Our two critical configurations of $H(T_{n,m})$ are perhaps a first step
toward proving this.

\subsection{Elastic tension energies}
All of the links presented above are critical or minimizing for the sum 
of the lengths of their components.  This is a beautiful functional, but
it is physically somewhat unrealistic: elastic ropes should minimize a 
quadratic functional of the form
\begin{equation*}
\sum_i a_i \big(\Len(L_i)-\ell_i\big)^2,
\end{equation*}
where $a_i>0$ is the elasticity and~$\ell_i$ the rest length
of the~$i$--th component.
Criticality for this functional is equivalent to that
for $\sum t_i \Len(L_i)$ where the tension~$t_i$
is $t_i=2a_i(\Len(L_i)-\ell_i)$.  Assuming these tensions
are nonnegative (that is, that no component's length is less than
its rest length) our balance criterion extends immediately
to handle this case: the strut force~$\Ag^*(\mu)$
must balance the tension-weighted curvature force $\sum t_i \K_i$.

In the known minimizing links, such as the simple chain,
each component separately achieves its minimal possible length.
Thus these examples also minimize all elastic energies with
nonnegative tensions $t_i\ge0$.

This behavior, however, seems rather exceptional.
The examples in \prop{wrapped} do not minimize all such functionals. 
In particular, if the tension in the long component is large 
enough, it will shrink to length $O(\sqrt{k})$ while some of the shorter
components gain length.  

Also in the Borromean rings, if the three components have different tensions,
the configuration we describe below (\secn{borromean})
would no longer be critical.
Similarly, clasps (\secn{gehr-clasp}) in which the two ropes
have different tensions again have new critical configurations.
In~\cite{SW} we describe in detail the shapes of these asymmetric clasp curves,
as well as their appearance in more complicated clasp-like links
even when tensions are equal.  (Note that \gehring thickness was
called Gehring thickness there, as in early drafts of this paper.)

\subsection{Nonembedded critical links}
To illustrate the differences between the Gehring problem
and the original ropelength problem, we now give some examples of
a different flavor: critical configurations that are nonembedded
and thus have infinite ropelength in the original sense.

Any knot is of course link-homotopic to the unknot.  The \gehring
ropelength minimizer degenerates to a point (of length zero).
The same happens for any component of an arbitrary link that is
link-homotopically split from the rest of the link.

Milnor showed that, up to link homotopy,
links of two components are classified by their linking number~\cite{MR17:70e}.
\begin{example}\label{ex:unstable}
When the linking number is zero, the components split,
and the \gehring ropelength minimizer degenerates to have length zero.
We can, however, also describe another critical configuration for
this unlink: one component degenerates to a point~$p$ while
the second is a unit circle centered at~$p$.  This is clearly
an unstable critical point: obvious deformations can decrease
the ropelength to second order.
\end{example}

\begin{figure}[ht!]\centering
\begin{overpic}[scale=1.25]{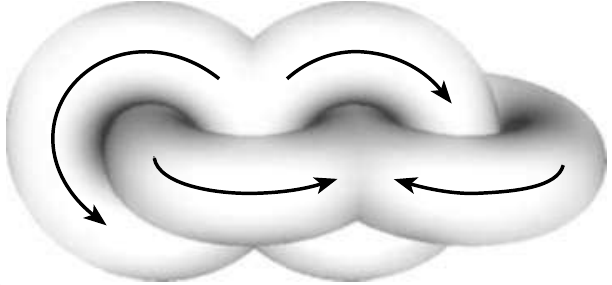}
\small
\put(21,40){$m_1$}
\put(56,40){$m_2$}
\put(38,11){$n_2$}
\put(78,11){$n_1$}
\end{overpic}
\caption[Another \gehring ropelength-critical link]
{In this configuration of two curves from \xmpl{covered-hopf},
each circle is covered~$m_i$ or~$n_i$ times,
as labeled. If $m = m_1 + m_2$ and $n = n_1 + n_2$, then
the curve has total length $2\pi(m+n)$, and linking number $mn-m_1n_1$.
It is constrained-critical, though often not minimal, for the Gehring
problem in the link-homotopy class defined by its linking number.}
\label{fig:m1n1}
\end{figure}

The case of linking number~$1$ is close to Gehring's original problem:
the minimizer is the same Hopf link built from round circles.  
(This case fits in the class $H(T_n)$ considered above.)
For larger linking number, we can use \cor{gehr-bal}
to exhibit many critical configurations as follows:

\begin{example}\label{ex:covered-hopf}
For linking number~$mn$ there is
a critical configuration~$L_{m,n}$ consisting of the minimizing
Hopf link with one component covered $m$~times and the
other covered $n$~times.  Its total length is thus $2\pi(m+n)$.
There are other critical configurations, sometimes shorter.
For example, each component can be a figure-eight built
from two tangent circles.
\figr{m1n1} shows a configuration like this with
total length $2\pi(m+n)$ and linking number $mn-m_1n_1$. 
The best configurations we know for linking number~$17$,
for instance, use $(m,n)=(6,3)$ or $(4,5)$.
Assuming configurations like these are the minimizers for
two-component links, they give examples where the set of
minimizers is disconnected (since we can interchange the two
components, or reorder the way one component covers its figure-eight).
\end{example}

\begin{figure}[ht!]\centering
\includegraphics[width=2.5in]{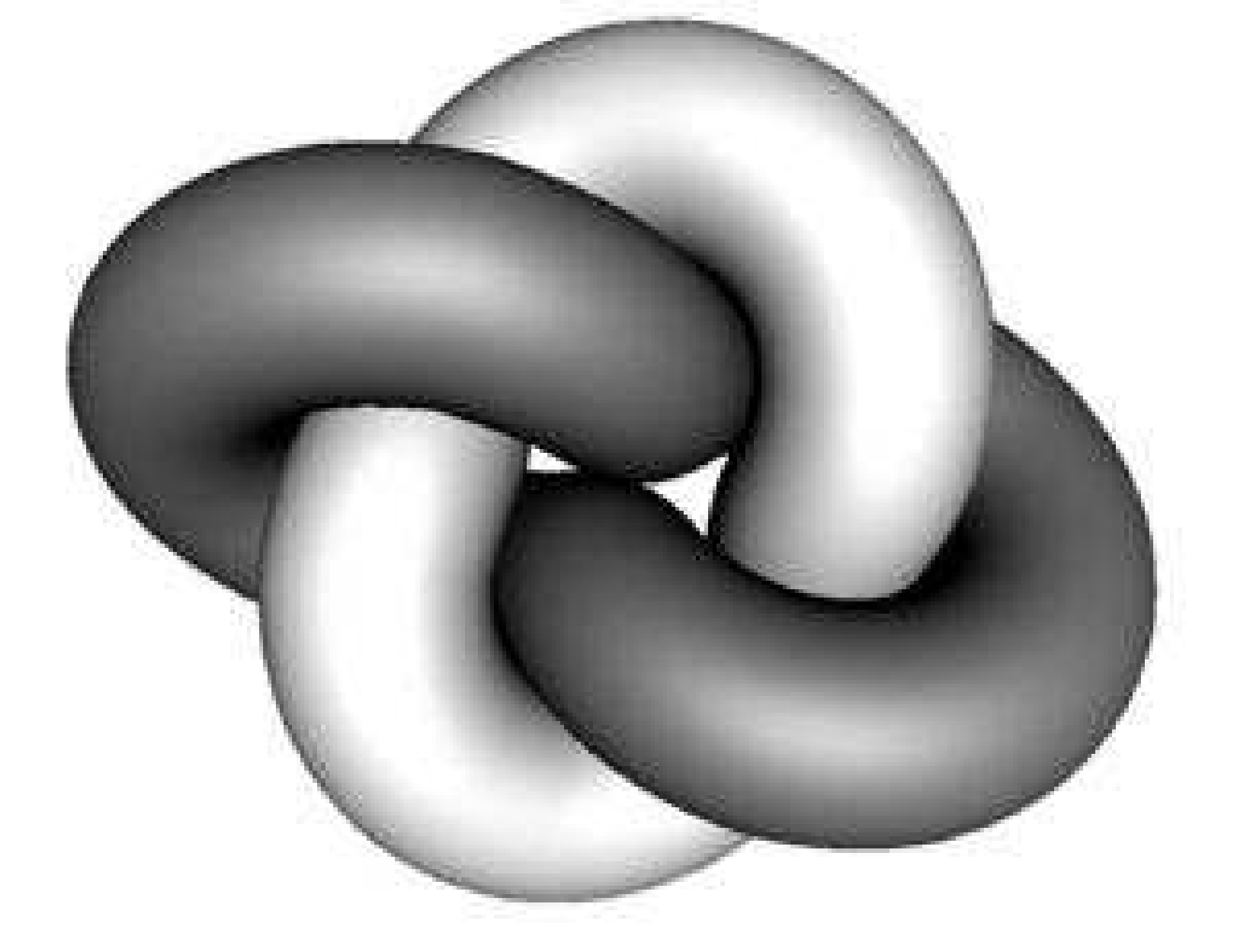}
\caption[The $(2,4)$--torus link]{This picture shows a numerically computed
minimizer for the \emph{original\/} ropelength problem on the $(2,4)$--torus
link. Because it has a strut between two points on the same component
(shown center, where the darker tube contacts itself) that carries
nonzero force, it is not
balanced for the Gehring problem considered here.
It is longer than~$L_{2,1}$, the Hopf link with one component doubly covered,
which we conjecture is the minimizer for \gehring ropelength
in this link-homotopy type.
Notice that both of these configurations break the symmetry between
the components of the link, so we expect to also find a (longer) critical
configuration where the two components are congruent.}
\label{fig:toruslink}
\end{figure}

None of these configurations is embedded, so they are not
critical points for the original ropelength problem:
as expected, the extra freedom in the Gehring problem sometimes allows
for shorter solutions.  As a further example, consider the $(2,4)$--torus
link, with linking number~$2$.  We have computed the presumed
ropelength-minimizer numerically, as in~\cite{ropen}. The results
are shown in \figr{toruslink}; this solution is longer than
the covered Hopf link $L_{2,1}$ (the presumed \gehring ropelength-minimizer)
and is not even critical for the Gehring problem.

For links of more than two components, linking numbers do not suffice
to distinguish link-homotopy types; we must also consider
Milnor's $\mu$--invariants \cite{MR17:70e,MR91e:57015}.
For instance, the Borromean rings, with no nonzero linking numbers,
belong to a nontrivial link-homotopy class because they have
$\mu$--invariant equal to~$1$.

Numerical experiments performed with Brakke's Evolver (compare \cite{ropen})
suggest that the minimizing Borromean rings for the \gehring ropelength
problem should consist of three congruent curves in perpendicular planes.
In~\cite{CKS2}, we described such a configuration built from circular
arcs of radius~$1$.  Unfortunately, \cor{gehr-bal} shows this configuration
is not even critical for length when constrained by \gehring thickness. 
In \secn{borromean}, the culmination of our paper, we will
explicitly describe a very similar configuration of the Borromean rings,
which we prove is critical and believe is the minimizer.

However, in order to solve for these Borromean rings, we must first
consider a simpler interaction between two ropes:
the clasp that occurs when one rope is pulled over another.
Describing this will require a notion of generalized links.

\section{Generalized link classes}\label{sec:genlink}

Although some of our definitions have applied to arbitrary
curves, so far we have been treating only ordinary (closed) links.
We now want to consider generalized problems involving curves
with endpoints.  To get meaningful link classes in this setting,
we must include constraints for the endpoints and obstacles for the link.

\begin{definition}
A \emph{generalized link\/}~$L$ is a curve~$L$
(with disjoint components $L_1,\dots, L_N$)
together with obstacles and endpoint constraints.
In particular,
each endpoint $x\in\bd L$ is constrained to stay on
some affine subspace $M_x\subset\R^3$, which can have
dimension~$0$, $1$ or~$2$.
Furthermore, there is a finite collection of \emph{obstacles\/}
for each component~$L_i$ of the link.
Each obstacle $$\{p\in\R^3: g_{ij}(p)<0\}$$ is given by a~$C^1$
function~$g_{ij}$ with~$0$ as a regular value.
By calling them obstacles, we mean that~$L_i$ is constrained
to stay in the region where $\min_j g_{ij} \ge 0$.
\end{definition}

While we could allow even more general endpoint and obstacle constraints,
this version fits nicely with our overall setup, and allows for all the
specific examples we have in mind.

\begin{definition} \label{defn:linkclass}
Suppose~$L= \bigcup_i L_i$ is a generalized link, with
obstacles~$g_{ij}$ and endpoint constraints~$M_x$.  Then its
\emph{link-homotopy class\/} $[[L]]$ is the set of all links~$L'$ that
are link-homotopic to~$L$ through links where each component avoids
its obstacles and maintains its endpoint constraints.  (As before, in
a link homotopy, each component of~$L$ can intersect itself but not
the others.)
\end{definition}

This definition is comparable to our previous definition for closed links
(\fullref{sec:gehringdef}); as in the discussion at the end of 
\secn{gehrbal}, we may restrict our attention to tame link classes.

Given a generalized link~$L$, only variations preserving
the endpoint constraints should be allowed.
A vector field $\xi\in\VF(L)$ is said to be \emph{compatible\/} with the
constraints if it is tangent to~$M_x$ at each endpoint $x\in\bd L$.
We write~$\VF_c(L)$ for the space of all compatible vector fields.

Given a set of obstacles $g_{ij}<0$ and a link~$L= \bigcup L_i$, we write
$$O(L) := \min_{i,j} \min_{x\in L_i} g_{ij}(x).$$
Then~$L$ avoids the obstacles~$g_{ij}$ if and only if $O(L)\ge 0$.
We define the wall struts of~$L$ by
$$\wall_{ij}(L) := L_i\cap\{g_{ij}=0\},\quad \wall(L):=\bigdisj_{i,j} \wall_{ij}(L).$$
This incorporates those parts of~$L$ on the boundary of the obstacle,
but is not strictly speaking a subset of~$L$ since one point $x\in L_i$
might be in several of the~$\wall_{ij}$.
When $O(L) = 0$, by Clarke's \thm{clarke} we have
$$\dplus_\xi O(L)
 = \min_{i,j}\min_{x\in \wall_{ij}(L)} \langle<\xi_x,\grad g_{ij}\rangle.$$
Again, we collect the various derivatives appearing on the right-hand side
into a rigidity operator $\Aw\co \VF_c(L)\to C(\wall(L))$ on wall struts,
given by
$$\Aw\xi(x):=\langle \xi_x,\grad g_{ij}\rangle$$
when $x\in \wall_{ij}$.
Its adjoint $\Aw^*$ is then
\begin{equation}\label{eq:Aw}
\int_L \xi\,\d \Aw^*(\mu)
  = \int_{\wall(L)}\!\! \Aw\xi\dmu
  =  \sum_{i,j} \int_{x\in \wall_{ij}(L)}\!
        \langle\xi_x,\grad g_{ij}\rangle\dmu(x).
 \end{equation}
We also have corresponding definitions for locally minimal,
strongly critical, and critical configurations of~$L$:

\begin{definition}
\label{defn:obstaclelocalmincritstrongcrit}
We say that a generalized link~$L$ is a \emph{local minimum\/}
for length when constrained by $\gThi$ if we have $\Len(L')\ge\Len(L)$,
for all sufficiently $C^0$--close links~$L'$ with the
same obstacle and endpoint constraints and with $\gThi(L')\ge\gThi(L)$.
We say~$L$ is \emph{strongly critical\/} (respectively, is \emph{critical\/})
for minimizing length when constrained by $\gThi$
if there is $\eps>0$ such that
for all compatible smooth~$\xi$ with $\delta_\xi \Len = -1$, the quantity
$$\min\big(\dplus_\xi \gThi(L), \,\dplus_\xi O(L) \big)$$
is at most~$-\eps$ (respectively, is negative).
\end{definition}

As in our discussion of Kuhn--Tucker
at the beginning of \secn{constrainedcritical},
these notions will be equivalent only under
a regularity assumption corresponding to the constraint qualification
of Mangasarian and Fromovitz~\cite{MR34:7263}:
\begin{definition}
\label{defn:obstaclemfcq}
A generalized link~$L$ is \emph{$\gThi$--regular\/} if there is
a \emph{thickening field\/}, meaning a smooth compatible~$\eta$
for which $\dplus_\eta \gThi(L)>0$ and $\dplus_\eta O(L)>0$.
\end{definition}

Note that, while we require~$\eta$ to strictly increase $\gThi$
and to move away from the obstacles, both to first order, there
is no corresponding requirement for the endpoint constraints,
since they are linear \emph{equality\/} constraints instead of nonlinear
inequality constraints.

We can now prove a generalization of \prop{gmincrit}:

\begin{proposition}\label{prop:th-regular}
If a generalized link~$L$ is a $\gThi$--regular
local minimum when constrained by $\gThi$, then~$L$ is critical. 
Also, if~$L$ is $\gThi$--regular and critical
when constrained by $\gThi$, then it is strongly critical.
\end{proposition}

\begin{proof}
The regularity of~$L$ means there exists a thickening field $\eta\in\VF_c(L)$.
We may assume $\delta_\eta\len(L)\ge0$ for otherwise~$L$ is neither minimal
nor critical; we then scale~$\eta$ so that $\delta_\eta \len(L) < 1$.

Suppose that~$L$ is a local minimum but not critical.
Then for some compatible vector field~$\xi$ we have
$\delta_\xi\len(L)<0$ while $\smash{\dplus_\xi}\gThi(L)\ge 0$
and $\smash{\dplus_\xi} O(L)\ge 0$.
For small $t>0$, consider the links $L_t=L+t(\xi+ \epsilon\eta)$.
Then 
\begin{equation*}
\ddtz{\len(L_t)}{+}{} = \delta_\xi \len(L) + \epsilon \delta_\eta \len(L).
\end{equation*}
We choose $0 < \epsilon < -\delta_\xi \len(L)/ \delta_\eta \len(L)$, so
this derivative is negative at time~$0$.  Thus for small~$t$, the~$L_t$
have length less than $\len(L)$, contradicting minimality if they obey
our constraints. But
\begin{equation*}
\ddtz{\gThi(L_t)}{+}{} \!> 0, \qquad 
\ddtz{O(L_t)}{+}{} \! > 0,
\end{equation*}
and the endpoint constraints are linear, so the links~$L_t$ meet all
our constraints for small $t>0$.

Now suppose that~$L$ is critical without being strongly critical.
Then there exists a sequence of compatible vector fields
$\xi_i \in \VF_c(L)$ with $\delta_{\xi_i} \len(L) = -1$ but with either
$\dplus_{\xi_i} \gThi(L) \approaches 0$ or
$\dplus_{\xi_i} O(L) \approaches 0$.
Then we observe that
$\delta_{\eta + \xi_i} \len(L) < 0$
for all~$i$, while by \cor{dplus-super} either
$$\lim\dplus_{\eta + \xi_i} \gThi \ge \dplus_{\eta} \gThi > 0$$
$$\lim\dplus_{\eta + \xi_i} O \ge \dplus_{\eta} O > 0.\leqno{\hbox{or}}$$
Taking~$i$ large enough that one of these quantities is positive,
we get a contradiction to the criticality of~$L$.
\end{proof}

So far, this development has paralleled that of \secn{gehrbal};
we now diverge from our previous course. Earlier, we saw that every
closed link is $\gThi$--regular:
rescaling always provides a thickening field. In the generalized setting,
this is no longer the case.
Thus minimality no longer implies criticality.
\begin{example}\label{ex:link-mnc}
To give a specific example,
rotate the constraints of \xmpl{min-not-crit} around the $z$--axis
to give obstacles $g_1=(x^2+y^2-1)^3-z<0$ and $g_2=z<0$ for an unknot~$L$.
The unit circle in the $xy$--plane is on the boundary of both obstacles,
and is clearly the minimum-length configuration in its homotopy class.
However, it is not critical: shrinking it toward the origin will
reduce its length to first order; the constraint $g_1\ge0$ is now violated,
but not to first order.
\end{example}

Further, criticality
and strong criticality may be different: if we allowed infinitely many
obstacles, we could construct critical, but not strongly critical links by
following the lead of \xmpl{strcrit-not-crit}.  (If we do
\emph{not\/} allow infinitely many obstacles, then an open question
remains: is strong criticality a stronger assertion than criticality?)


\begin{example}\label{ex:link-sc-not-reg}
To justify our emphasis on strong criticality (rather than restricting
our attention to regular, critical links) we also note that it is easy
to construct strongly critical links that are not regular; simply
take~$L$ to be the unit circle in the $xy$--plane, with constraints
$g_1(x,y,z) = x^2 + y^2 - 1$ (so the excluded region is the infinite
cylinder around the $z$--axis) and $g_2 = -g_1$.
This link is trapped on the cylinder $g_1=0=g_2$,
so it has no thickening field. On the other hand, it is clearly
strongly critical.
\end{example}

Now we are ready to extend our balance theorem
to the generalized setting.
We will accommodate the endpoint constraints by restricting our
attention to compatible vector fields.  Our other constraints are then
$\dist\ge1$ on~$\UPL$ and $g_{ij}\ge 0$ along~$L_i$.
The set~$Y$ of active constraints then consists
of the struts together with the wall struts.

\begin{theorem}\label{thm:gen-gehr-bal}
A generalized link~$L$ is strongly critical for \gehring ropelength
if and only if there is a positive Radon measure~$\mu$ on
$\gStrut(L)\disj\wall(L)$, such that $$-\K=(\Ag \oplus \Aw)^*\mu$$ 
as linear functionals on $\VF_c(L)$.  This means
that~$-\K$ and $(\Ag \oplus \Aw)^*\mu$ agree as forces along~$L$
except at endpoints $x \in\bd L$, where they may differ by an 
atomic force in a direction normal to~$M_x$.
\end{theorem}

\begin{proof}
This is again a straightforward application of our \thm{ourkt}, using
$$X=\VF_c(L), \quad Y=\gStrut(L)\disj\wall(L), \quad
f=\delta\,\len, \quad A=\Ag\oplus\Aw .\proved$$
\end{proof}

\begin{remark}
Remember that $\K$ has been defined to include an inward-pointing
atom at each endpoint $x\in\bd L$.
We can ignore these, however, when applying this theorem,
as long as the link~$L$ meets each endpoint constraint~$M_x$ normally.
We know of no examples of critical links where this is not the case.
\end{remark}

The regularity described in \cor{gehr-ftc} and \prop{gcrit-c1}
carries over to generalized links:

\begin{proposition} \label{prop:gen-gcrit-c1}
If the generalized link~$L$ is strongly critical for \gehring ropelength,
then~$L$ is~$\FTC$ and~$C^1$.
\end{proposition}

\begin{proof} 
The proof follows that of \cor{gehr-ftc} and \prop{gcrit-c1}.
From equation~\eqref{eq:Aw}, we find that $\Aw^*\mu$
has distributional order zero just like $\Ag^*\mu$,
so $L\in\FTC$ follows immediately from the balance criterion
of \thm{gen-gehr-bal}.

Now suppose $L$ is not~$C^1$ but instead has some corner~$x$
with $\K\{x\}\ne 0$.  By \mbox{\thm{gen-gehr-bal}}, this curvature
force is balanced by struts and wall struts.  So there
is at least one strut or wall strut acting on~$x$
in a direction with negative inner product with $\K\{x\}$.
In the case of a strut $\{x,y\}$, we refer again to \figr{gehr-corn}:
some points near~$x$ along~$L$ would be nearer to the endpoint~$y$.
But similarly, in the case of a wall strut, we have
$\langle\K,\grad g_{ij}\rangle< 0$, but this means that
some points near~$x$ along~$L$ violate this obstacle constraint. 
In either case, we get the desired contradiction.
\end{proof}

To understand the interplay between struts and wall struts, we now
offer a simple example of a generalized link~$L$ with
nonempty boundary which is balanced, needing nonzero force on the wall struts.
\begin{example}\label{ex:wall}
Cut the simple chain of \fullref{fig:chainstruts} by parallel
planes through~$s_1$ and~$s_2$ with normal vector $c_1-c_2$, and let~$L$
be the part of the chain lying between the two planes.
This generalized link includes two semicircles with endpoints normal
to the planes, and also the inner stadium curve, which is tangent
to the planes at~$s_1$ and~$s_2$.  We let the planes bound an obstacle,
forcing~$L$ to stay between the planes, and we use them also as endpoint
constraints.  Then~$L$ is balanced: though the semicircles
now exert a net outward force on~$s_1$ and~$s_2$, this is balanced by 
wall struts at these points. And the internal balance for the stadium
curve remains the same.
\end{example}

\section{The tight clasp}\label{sec:gehr-clasp}

The tight configurations of \secn{gehr-examples} were the
simplest closed links we could imagine:
the Hopf link, and various connect sums of Hopf links in which
each component is still a convex plane curve. 
But there is an even simpler interaction between
two ropes, the \emph{clasp\/} formed when one rope is pulled taut over
another, as at the junctions of a woven net, or when a bucket is lifted
from a well by passing a rope through its rope handle.
We can model a single clasp as a generalized link with endpoint constraints.

To define the \emph{simple clasp\/}, fix two parallel planes~$P$
and~$\tild P$ at least~$2$ units apart.
Then take two unknotted arcs~$\gamma$ and~$\tg$ that
lie between the planes, with the endpoints of~$\gamma$ constrained to lie
in~$P$ and those of~$\tg$ in~$\tild P$.
Let the halfspace bounded by $P$ that does not include~$\tild P$ be an obstacle
for the component $\tg$, and vice versa, and
select the isotopy class of such links shown in \figr{clasp}.
This is the class where
closing each arc in the plane of its endpoints would produce a Hopf link.

\figs{clasp}{1.25}{\figdir/looseclasp}{The simple clasp}
{The simple clasp has two components, one attached at both ends to
the ceiling and the other to the floor, linked with one another as shown.
The configuration shown, with each component consisting of a semicircle joined
to two straight segments, is neither critical nor minimal.}

It is natural to assume that the minimizing configuration for
this problem would consist of semicircular arcs passing through
each others' centers, together with straight segments joining the
semicircles to the constraint planes, much like the Hopf chain of
\xmpl{hopf-chain}.  But this \emph{naive clasp\/} is not balanced:
each semicircle focuses its curvature force on the tip of the other,
and there is no way to balance these forces (as the isolated strut
carrying an atom of compressive force did in the Hopf chain).
The naive clasp is thus not minimizing, though we will see it is very close:
the critical configuration we construct here is only half a percent shorter.

\begin{example}\label{ex:pressed-clasp}
Suppose the horizontal planes~$P$ and~$\tild P$ in the definition of the
simple clasp are taken instead to be only one unit apart.  Consider the
configuration where the curves~$\gamma$ and~$\tg$ are semicircles
in perpendicular vertical planes.
The curvature of each semicircle can be balanced by uniform strut tension,
transmitting a net vertical force to the tip of the other semicircle.
That vertical force can be balanced by a wall strut at each tip.
Therefore, this configuration is critical for \gehring ropelength.
\end{example}

In the case of interest, where~$P$ and~$\tild P$ are far apart and there
are no wall struts to balance the tips,
we must look harder for a solution.
We will now construct critical configurations,
constrained by the \gehring thickness $\gThi$, for the simple clasp problem
and for a family of related problems where the ends of the ropes are
pulled outward as in \figr{angleclasp}.
These solutions minimize length under natural symmetry assumptions, and we
believe they are the global minimizers even without imposed symmetry.
Below in \secn{borromean}, we will construct a critical configuration
of the Borromean rings that contains portions of these clasp curves.
Thus, a thorough understanding of these generalized links will aid us
in understanding that more complicated closed link.

\subsection{Symmetry conditions and a convenient parametrization}
We describe configurations of the clasp where the two components
are congruent plane curves, lying in planes perpendicular to
each other and to the constraint planes.  To fix these symmetries
in coordinates, let the constraints be the planes $z = \pm C$, and let the
component~$\gamma$ lie in the $xz$--plane while~$\tg$ lies in the $yz$--plane.
The clasp has mirror symmetry across each of these planes (preserving each
component).  It also has a symmetry interchanging the two components,
which we denote $p\mapsto\tild p$, given by
fourfold rotation about the $z$--axis together with reflection across
the $xy$--plane.  These symmetries generate a point group of order eight
in~$O(3)$ whose Conway--Thurston orbifold notation
\cite{conway3d,conway48} is~\nstartwo{2}.
Algebraically it is isomorphic to~$D_4$.

The argument we present below to derive the critical
clasps for the Gehring problem can easily be extended to show these are the
unique critical configurations among curves with this~\nstartwo{2} symmetry.
We omit the details, however, because we know of no way to show that
the overall minimizers must have this symmetry. If one could prove this,
it would then follow that our clasps are the minimizers.

Our symmetry assumptions mean that the clasp is described by the shape
of half of the component~$\gamma$, from its \emph{tip\/} along the $z$--axis
into the $x > 0$ half-plane and up to the plane~$P$.  This consists of
a curved arc near the tip joined to a straight segment near~$P$.
Since the curved arc is strictly convex, we can parametrize it
by the angle~$\phi$ made by
its tangent vector above the horizontal, as in \figr{orthoclasp}.
In fact, we will use the sine of this angle, $u=\sin\phi$, as our parameter.
Thus in the simple clasp, for $u \in [0,1]$ we write
\begin{align*}
\gamma(\pm u) &= (\pm x(u), 0, z(u)), \\
\tg(\pm u) &= (0, \pm x(u), -z(u)).
\end{align*}
Elementary calculations show the following:
\begin{lemma}\label{lem:cvx-crv}
For a convex curve~$\gamma$ in the $xz$--plane, parameterized by the
sine~$u$ of its direction~$\phi\in[-\tfrac\pi2,\tfrac\pi2]$,
the arclength~$s$ satisfies
$$ \d s = \sec\phi \dx = \csc\phi \dz = \frac{\d u}{\kappa \sqrt{1 - u^2}},$$
where the curvature~$\kappa$ is given by
$$\kappa=\frac{\d\phi}{\d s} = \frac{\d u}{\d x}.$$
\end{lemma}

For the simple clasp described above, each component turns a total
of~$180^\circ$, meaning that~$u$ ranges from~$-1$, through~$0$ at the tip,
to~$1$.  We can also consider more general clasp problems where the four
ends of rope are not vertical (being attached to horizontal planes)
but instead are pulled out at some angle (being attached to tilted planes).  

\begin{figure}[ht!]\centering
\begin{overpic}[scale=1.25]{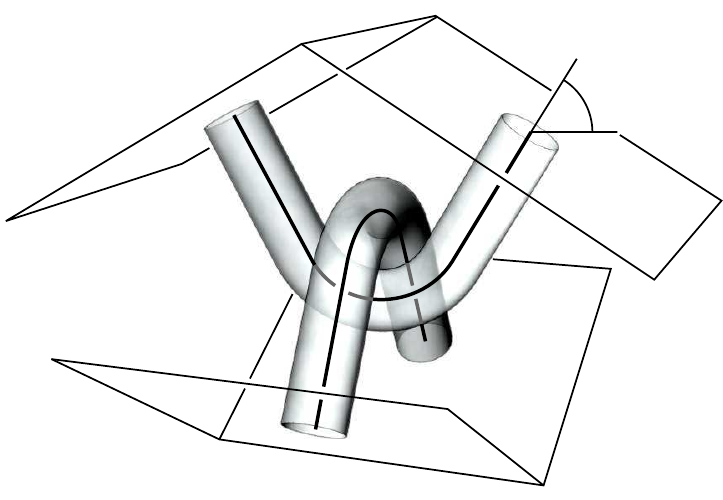}
\small
\put(83,56){$\arcsin\tau$}
\end{overpic}
\caption[The $\tau$--clasp problem] {In this variant of the simple
clasp problem, the endpoints of the two ropes are constrained to lie
in four planes whose normals make angle $\arcsin\tau$ with the
horizontal.  The parameter $u=\sin\phi$ ranges from~$-\tau$ to~$\tau$
along each arc, as shown at the end of the top right arc.
If extended, the four planes shown would form the sides of a tetrahedron.
Each arc is constrained to lie in the wedge formed by the planes containing
the endpoints of the other arc.}
\label{fig:angleclasp}
\end{figure}

Given $0\le\tau\le1$, we define the $\tau$--clasp to be a problem like
the simple clasp where the arc~$\gamma$ starts at $u=-\tau$ and then
turns through angle $2\arcsin\tau$ to reach $u=\tau$.  Our critical
$\tau$--clasps have the same \nstartwo{2} symmetry as the simple clasp.
To put the $\tau$--clasp into our framework of generalized links, we
constrain the four endpoints to four planes, each making angle
$\arcsin\tau$ with the vertical, as in \figr{angleclasp}. The
complement of the wedge formed by the planes containing the endpoints
of each arc acts as an obstacle for the other arc.  The simple clasp
is the $\tau$--clasp with $\tau=1$, where the wedges degenerate to
halfspaces.

\subsection{Struts between perpendicular planes}
Whenever two curves in perpendicular planes are connected
by a strut, elementary trigonometry gives us first order
information about the curves at both endpoints.  We state
a general lemma, which we will use here for the clasp
and again for the Borromean rings.

Let~$P_1$ and~$P_2$ be two planes meeting perpendicularly along
a line~$\ell$, and let $\gamma_i\subset P_i$ be two components of a link.
At a point $p_i\in\gamma_i$, we write~$x_i$ for the distance
from~$p_i$ to~$\ell$, and~$u_i$ for the cosine of the angle
between~$\ell$ and the line tangent to~$\gamma_i$ at~$p_i$.  These quantities
generalize the~$x$ and~$u$ of \lem{cvx-crv} above.

\begin{lemma}\label{lem:perp-planes}
Let~$\gamma_1$ and~$\gamma_2$ be two components of a link~$L$,
lying in perpendicular planes.  Suppose there is a strut $\{p_1,p_2\}$
of length~$1$ connecting these components.  Then in the notation
of the previous paragraph we have $0\le x_i\le u_i\le 1$,
and any two of the numbers $x_1,x_2,u_1,u_2$ determine the other two,
according to the formulas
\begin{align*}
x_i^2 &= 1-\frac{x_j^2}{u_j^2} = \frac{u_i^2(1-u_j^2)}{1-u_i^2u_j^2}, \\
u_i^2 &= \frac{1-x_j^2/u_j^2}{1-x_j^2} = \frac{x_i^2}{1-x_j^2},
\end{align*}
where $j\ne i$.
\end{lemma}
\begin{proof}
Picking cartesian coordinates such that~$\ell$ is the $z$--axis
and~$P_i$ are coordinate planes,
we find the strut difference vector $p_1-p_2$ is $(x_1,x_2,\Delta z)$,
for some number~$\Delta z$.
Since this strut has length~$1$ and is perpendicular to each~$\gamma_i$,
we have
$$ \Delta z^2 + x_1^2 + x_2^2 = 1,
\quad \Delta z = x_i \frac{u_i}{\sqrt{\smash[b]{1-u_i^2}}}. $$
Simple algebraic manipulations, eliminating~$\Delta z$,
lead to the equations given.
\end{proof}

Note that the condition $x_i\le u_i$ is exactly the condition that
the unit normal circle around~$p_i$ intersects~$P_j$; the two
points of intersection are mirror images (across~$P_i$),
with the same~$x_j$ and~$u_j$ values.
Also note that we don't need to have $\gamma_i\subset P_i$ in the lemma;
it suffices that~$\gamma_i$ be tangent to~$P_i$ at~$p_i$.

Whenever we have a pair of curves in perpendicular planes, which
stay a constant distance~$1$ apart, we can apply this lemma everywhere
along the curves.  Each curve~$\gamma_i$ is determined as the intersection
of the plane~$P_i$ with the unit-radius tube around the other curve~$\gamma_j$.
This will be the situation for the clasp.

\subsection{The balancing equations for the clasp}
By \thm{gen-gehr-bal}, in a critical clasp the curvature
force along~$\gamma$ must be balanced by struts to~$\tg$.
In particular, almost every point
(indeed, since the set of struts is closed, \emph{every\/} point)
$\gamma(u)$ along the curved arc of~$\gamma$
must have a strut to some point $\tg(u^*)$.
Then by symmetry we actually have what we call $2$--to--$2$ contact:
there are struts from~$\gamma(\pm u)$ to~$\tg(\pm u^*)$.
Here the two points $\tg(\pm u^*)$ must be the intersection
of the unit normal circle around $\gamma(u)$ with the $yz$--plane,
implying that $u^*\in[0,1]$ is uniquely determined for each~$u$.
We will refer to $\gamma(u)$ and $\tg(u^*)$
as \emph{conjugate points\/} on the $\tau$--clasp.
\lem{perp-planes} applies to any pair of conjugate points,
with $u_1=u$, $u_2=u^*$ and $x_i=x(u_i)$.

\begin{lemma}\label{lem:vert-force}
Suppose~$\gamma$ is a plane curve, symmetric across a line~$\ell$
in the plane.  Consider the net curvature force of a mirror image pair of
infinitesimal arcs of~$\gamma$.  This acts in the direction of the line~$\ell$,
with magnitude $2|\d u|$.  Here the function~$u$ is defined along~$\gamma$
as the cosine of the angle~$\psi$ between~$\ell$
and the tangent line to~$\gamma$.
\end{lemma}
\begin{proof}
One infinitesimal arc has net curvature force
$\kappa N\ds = N\,\d\psi$.  When this is added to the mirror
image force, only the component along~$\ell$ survives.
We get magnitude $2|\sin\psi\,\d\psi|=2|\d u|$.
\end{proof}

Suppose now we have a symmetric configuration of the clasp
where the curved arcs of the two components
stay a constant distance~$1$ apart.
By symmetry we get the $2$--to--$2$ strut pattern described above.
Assuming the straight ends of each component meet the constraint
planes perpendicularly, our balance criterion \thm{gen-gehr-bal}
says that strong criticality is equivalent
to the statement that the net vertical curvature force
exerted by the arcs at $\gamma(\pm u)$
balances that of the conjugate arcs at $\tg(\pm u^*)$.
That is, using \lem{vert-force}, for a critical clasp we must have
$|\d u|=|\d u^*|$, meaning that either $u-u^*$ or $u+u^*$ is constant.

If $u-u^*$ were constant, by symmetry it would be zero,
and our equations would describe
a pair of half-ellipses, with horizontal major axis~$\sqrt{2}$ and vertical
minor axis~$1$.  On these curves, corresponding points $\gamma(u)$ and $\tg(u)$
are always at distance~$1$ from each other, but these are \emph{maxima\/} for the
distance between components, rather than minima.
This configuration has $\gThi<1$, and is not $\gThi$--critical:
the pairs $\{\gamma(u),\tg(u)\}$ are not struts.

Instead we must have that $u+u^*$ is constant.
To find the constant, note that on the $\tau$--clasp,
the tip of~$\gamma$ (at $u=0$) is joined by a strut
to the end of~$\tg$ (at $u^*=\tau$); thus $u+u^*=\tau$.
This equation holds when $0\le u,u^* \le \tau$; to allow for
negative values (parametrizing the whole clasp curve) we write
$$|u|+|u^*|=\tau.$$
We can now give an explicit description of our critical $\tau$--clasp:

\begin{theorem}
\label{thm:clasp}
Let $\tau \in [0,1]$, and let $\gamma= \gamma_\tau$ be the curve in the
$xz$--plane given parametrically for $u\in[-\tau,\tau]$ by
\begin{align*}
x = x_\tau(u) &:= \frac{u \sqrt{1- (\tau-|u|)^2}}{\sqrt{1-u^2(\tau-|u|)^2}} , \\
z = z_\tau(u) &:= \int\frac{\d z}{\d x}\dx 
  = \int \frac{u}{\sqrt{1-u^2}} \,\frac{\d u}{\kappa_\tau(u)},\\
\tag*{\hbox{where}} 
\kappa_\tau(u) :=& 
   \frac {\sqrt{\big(1-u^2(\tau-|u|)^2\big)^3 \big(1-(\tau-|u|)^2\big)}} 
         {1-(\tau-|u|)^2+(\tau-|u|)|u|(1-u^2)}
\end{align*}
and the constant of integration for~$z$ is chosen so that
$$z(0)+z(\tau)=-\sqrt{1-\tau^2}.$$
Then the union of~$\gamma$ with its image~$\tg$ under the symmetry group
\nstartwo{2} described above is a $\tau$--clasp that is strongly critical
for \gehring ropelength.
The curvature of~$\gamma$ is $\kappa_\tau(u)$ above, and the
total length of the curved part of~$\gamma$ is
$$\int_{-\tau}^{\tau} \frac{\du}{\kappa_\tau(u) \sqrt{1 - u^2}}.$$
\end{theorem}
\begin{proof}
The proposition follows from the foregoing discussion, after
substituting $u^*=\tau-|u|$ into the equations of \lem{perp-planes},
and using \lem{cvx-crv}.  To get the constant of
integration for~$z$, we note that the strut from $\gamma(0)$ to $\tg(\tau)$
has height given (as in the proof of \lem{perp-planes}) by
$$\Delta z=\tsty\sqrt{\smash[b]{1-x_\tau(0)^2-x_\tau(\tau)^2}}=\sqrt{1-0-\tau^2}.\proved$$
\end{proof}

\begin{figure}[ht!]\centering
\begin{overpic}[width=.76\textwidth]{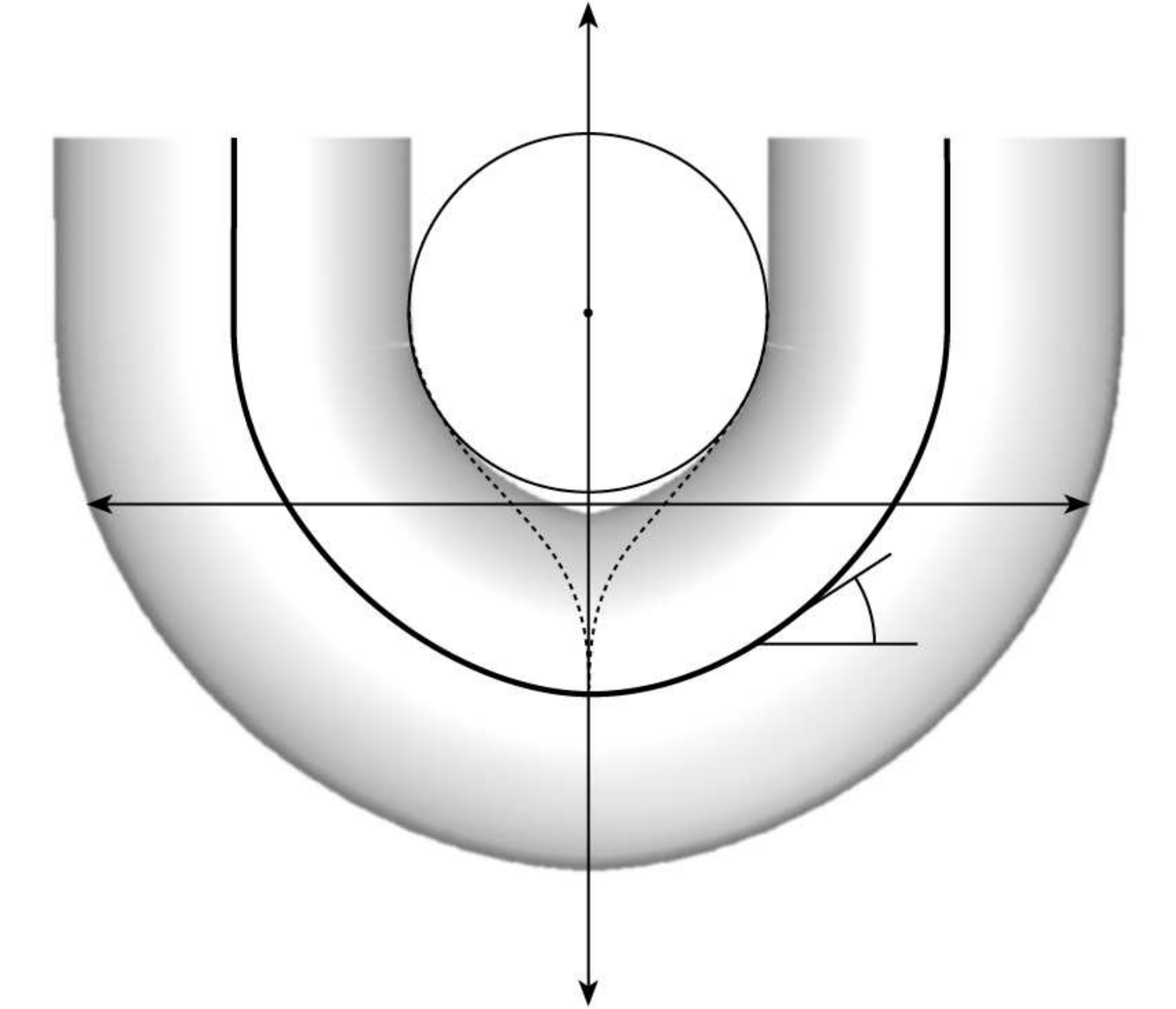}
\small
\put(85,44){$x$}
\put(46,80){$z$}
\put(75,34){$\phi$}
\put(51,22){$\gamma(0)$}
\put(51,58){$\tg(0)$}
\put(82,58){$\gamma(1)$}
\end{overpic}
\caption[The tight clasp]
{This is an accurate plot of the critical simple clasp~$\gamma$ given by
\thm{clasp}.  Here $u = \sin \phi$ ranges from~$-1$ to~$1$ over the curved
portion of~$\gamma$. The tip $\tg(0)$ of the other component is shown
above~$\gamma$ on the $z$--axis, along with the (dotted) circular cross-section
of the tube of unit diameter around~$\tg$. The curved dotted lines extending
down from the sides of this cross-section are the lines of contact between the
shaded tube around~$\gamma$ and the front half of the tube around~$\tg$.
Symmetric lines
of contact extend behind the shaded tube, realizing the $2$--to--$2$ contact
pattern we have described.  Finally, we see a small gap between the
tubes, explored in more detail in \figr{gapchamber}.}
\label{fig:orthoclasp}
\end{figure}

Although the formulas we have given for $z_\tau(u)$ and for arclength
both involve hyperelliptic integrals not expressible in closed form,
it is straightforward to integrate them numerically;
we have plotted our critical configuration
of the simple ($\tau=1$) clasp in \figr{orthoclasp}.

As we mentioned in the introduction, Starostin has
given~\cite{starostin-forma} an independent derivation (using
a form of balancing for smooth curves) of these same $\tau$--clasp
configurations (as well as the family of stiff clasps we will consider
in~\cite{CFKSW2}). Starostin does not prove that these
configurations are critical for \gehring ropelength.

\subsection{The geometry of the tight clasp}
We now examine the curvature
and other geometric features of the critical clasps for the
Gehring problem that were given in \thm{clasp}.
Each component of the critical $\tau$--clasp is a~$C^1$ join of four
analytic pieces: a straight segment, then $\gamma[-\tau,0]$,
then $\gamma[0,\tau]$, and finally another straight segment.
Where the curved arcs join the straight segments at $u=\pm\tau$,
the curvature~$\kappa(u)$ approaches~$1$;
at these points, our critical clasp agrees to second order
with the naively expected circular arcs.

\begin{figure}[ht!]\centering
\begin{overpic}[scale=0.5]{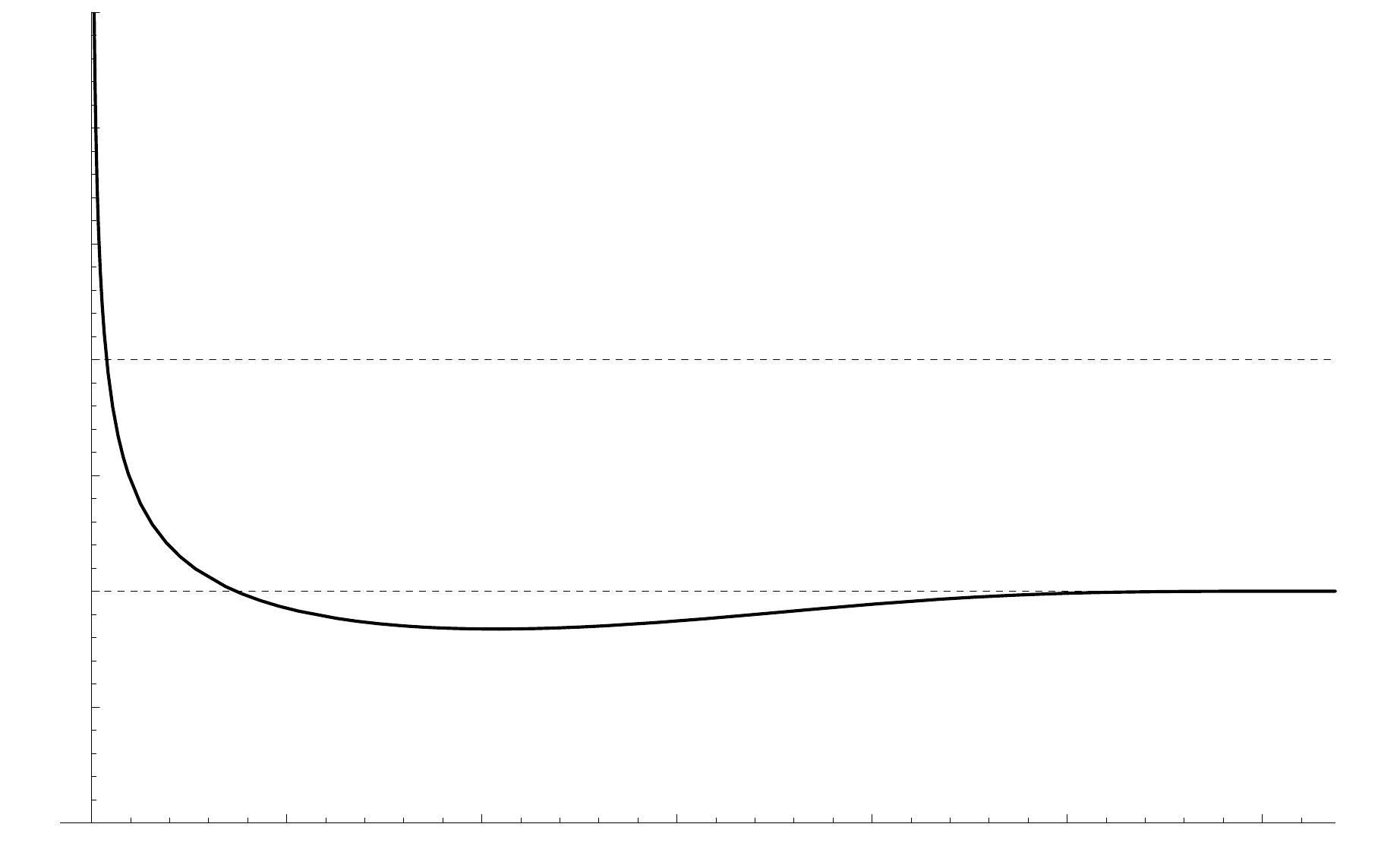}
\small
\put(50,-3){$s$}
\put(32,-1){$0.5$}
\put(61.5,-1){$1$}
\put(88,-1){$1.5$}
\put(0,30){$\kappa$}
\put(3,18){$1$}
\put(3,35){$2$}
\put(3,52){$3$}
\end{overpic}
\caption[The curvature of the tight clasp]
{The graph shows the curvature~$\kappa$ of the tight $1$--clasp as
a function of arclength. The curvature blows up at the tip: this curve is
only $C^{1,\nicefrac{2}{3}}$. The unit-diameter thick tube around
the curve forms a cusp near the tip, when the curvature exceeds~$2$.
From the tip, curvature decreases rapidly to its minimum,
and then increases again to the limiting
value of $\kappa=1$ at the end.  Thus the clasp curve,
at its end, agrees to second order with the naively expected unit circle
around the tip of the other component, as is suggested in \figr{orthoclasp}.
(For $\tau=1$, as illustrated here, the curves agree even to third order.)}
\label{fig:kappagraph}
\end{figure}

The maximum curvature $\kappa(0)=1/\sqrt{1-\tau^2}$ occurs at the tip.
For $\tau<1$, this is finite, and our $\tau$--clasp is $C^{1,1}$.
But for $\tau=1$, the curvature blows up (like $|s|^{-\nicefrac{1}{3}}$) at the tip.
In \figr{kappagraph} we plot the curvature $\kappa(u)$ for this
simple clasp.  The curve is $C^{1,\nicefrac{2}{3}}$ (and is also in the Sobolev
space $W^{2,\,3-\eps}$ for all $\eps>0$) but has no higher regularity.

\figthree{gapchamber}{.95}
{\figdir/tightclasp_expl2}{\figdir/tightclasp_transparent}{\figdir/tightclasp_closeup}
{The gap chamber in the clasp}{We see three views of the gap chamber
between the two tubes in the tight clasp with $\tau = 1$. On the left,
we see an exploded view with the two tubes and the gap chamber floating
between them. In the medium closeup in the center, we see the chamber in place
between the (now transparent) tubes. On the right, we see an extreme closeup
of the center of the gap chamber. Its height at the center (about~$0.05639$)
is the distance between the tubes at the tips of the clasp. The grid
in the center and right pictures is a square grid projected from the
$xy$--plane. On the right, we see a tiny ridge running from left to right
along the surface of the chamber; this is a cusp formed by the folding
of the tube surface that happens when the curvature of the clasp rises
above~$2$ (compare \figr{kappagraph}).  We do not know whether this gap chamber
forms in clasps of physical rope; it would be very interesting to find out.}

In \prop{gcrit-c1}, we proved that critical curves for \gehring ropelength
are~$C^1$.  It would be interesting to find out whether all such critical
curves are $C^{1,\nicefrac{2}{3}}$; perhaps the simple clasp exhibits the worst
possible behavior.

In \xmpl{covered-hopf}, we saw critical curves constrained
by \gehring thickness which
fail to have positive thickness in the ordinary sense of~\cite{CKS2}
because one component is nonembedded.  The simple clasp fails 
to have positive thickness for a different reason: its curvature is
unbounded.  In~\cite{CFKSW2} we will consider a family of thickness
measures with a variable stiffness parameter~$\lambda$. In these measures,
a unit-thickness curve has curvature bounded above by~$\nicefrac{2}{\lambda}$.
For any nonzero~$\lambda$, it follows that the critical simple clasp must be
different from the tight clasp here for the Gehring problem,
and must instead include an arc of this maximum allowed curvature.

One of the most interesting features of the clasp is the gap between
the two components of the clasp. The distance between the tips of~$\gamma$
and~$\tg$ is $z(\tau)-z(0)+\sqrt{1-\tau^2}$ (written in this way to
be independent of the constant of integration for~$z$). This is an increasing
function of~$\tau$, close to~$1$ when~$\tau$ is small, but increasing
to~$1.05639$ at $\tau=1$. Thus, in the simple clasp, the gap
between the thick tubes around the two components at their tips is
almost~$6\%$ of their diameter.

\begin{figure}[ht!]
\centering
\begin{overpic}[clip,scale=0.6]{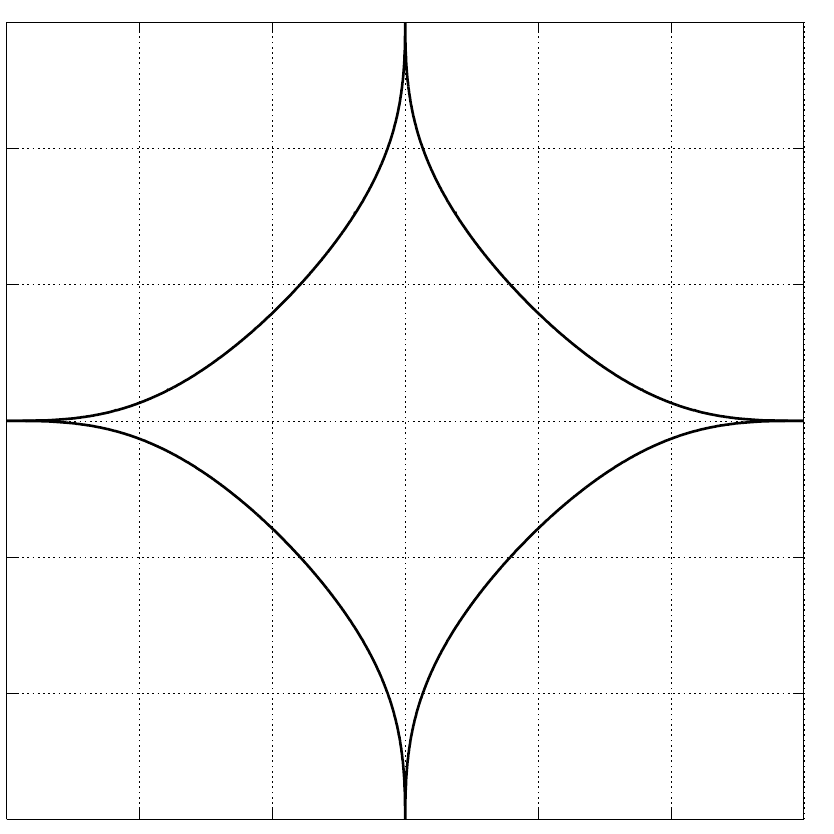}
\small
\put(-6,-5){$-\ell$}
\put(46,-5){$0$}
\put(96,-5){$\ell$}
\end{overpic}
\caption[The strut set of the tight clasp]
{The graph shows the strut set for the tight $1$--clasp, where each
strut is plotted according to the arclength of its ends
on the two components of the clasp (measured from the tip $s=0$
to the shoulders at $s=\ell\approx1.58944$).
There is a closed loop of struts,
with four cusps at the tips and shoulders of the clasp arcs.
We hope this explicit strut set will help in verifying
the accuracy of numerically computed strut sets for ropelength minimizers,
such as those in~\cite{mglob}.}
\label{fig:claspstruts}
\end{figure}

These thick tubes contact each other at the midpoints of the struts.
Topologically, the set of struts forms a loop.  Their midpoints
form a loop in space with four vertical cusps---%
the line of contact of the two tubes---%
as seen in \fullref{fig:orthoclasp} and~\fullref{fig:gapchamber}.
Alternatively, we can plot the loop
of struts as pairs of arclength coordinates on the two components,
as in \fullref{fig:claspstruts}.
The solid tubes divide the rest of the ambient space into two regions:
one infinite component around the outside of the clasp,
and one small chamber sitting in the gap between the tips,
shown in \figr{gapchamber}.  To give a sense of scale, the gap chamber has a
substantial surface area of about~$1.10$, equal to the area of a section
of tube of length more than~$\nicefrac{1}{3}$. However, the chamber is very thin,
resulting in a volume of only~$0.01425$.

\subsection{Length comparison with the naive clasp}
Earlier, we described the naive circular configuration for the simple clasp.
Similarly, in what we call the \emph{naive $\tau$--clasp\/}, each component
is built from straight segments (normal to the constraint planes) and a
unit-radius arc (of angle $2\arcsin\tau$ and centered
at the tip of the other component).
As we saw for $\tau=1$, this configuration is not critical:
there is no way to balance the forces concentrated on the tips,
unlike in \fullref{ex:hopf-chain} and~\fullref{ex:pressed-clasp},
which had extra struts.

Our critical $\tau$--clasps (which we expect are the global minima for length)
are indeed slightly shorter than the naive configurations.
The total length of a clasp depends, of course,
on the position of the bounding planes.
Thus to compare the lengths of the naive clasp and our critical clasp
in a meaningful way, we introduce the notion of excess length.
The infimal possible length
of a $\tau$--clasp with no thickness constraint is easily seen to be
four times the inradius of the bounding tetrahedron.  (In the case
$\tau=1$ this is twice the thickness of the bounding slab.)
The \emph{excess length\/} of any given clasp is the amount by which its
length exceeds this value.

For $\tau=1$, the naive clasp has excess length $2\pi-2$, since two
unit semicircles replace two straight segments of unit length.
Numerical integration reveals the excess length of our critical
$1$--clasp to be~$4.262897$ (accurate to the number of digits shown).
It is thus about~$0.020288$, or almost half a percent, shorter.
In general, the excess length of the naive $\tau$--clasp is
$4\arcsin\tau - 2\tau$, while the excess length of our critical
$\tau$--clasp equals the total length of the curved parts minus~$2\tau$
times the inter-tip distance.  The maximum percentage savings,
about~$0.518\%$, occurs for $\tau\approx\sin(80^\circ)$.

\section{The Borromean rings}\label{sec:borromean}

The original Gehring link problem was solved by the Hopf link made from a 
pair of circles through each other's centers. We have already generalized
this to a three component link in one way: the simple chain made 
from circles and stadium curves of \fullref{sec:gehr-examples}.
But the simple chain is just a connect sum of Hopf links, and so
the minimizing configuration shares much of its geometry with
the original Gehring solution.  

\begin{figure}[ht!]
\begin{center}
\includegraphics[scale=.75]{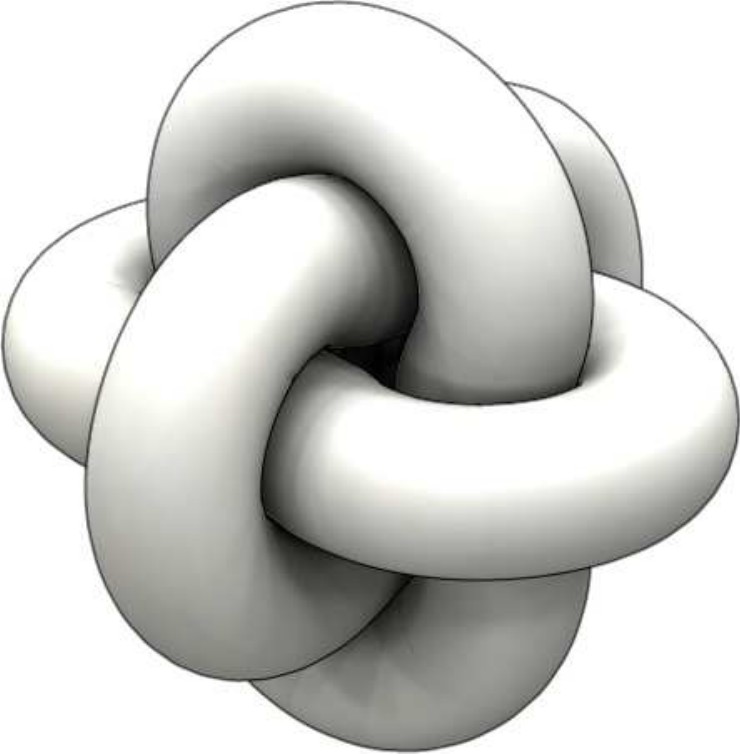}
\end{center}
\caption[The Borromean Rings]
{Our critical configuration~$B_0$ of the Borromean rings,
shown with thick tubes of diameter~$1$.
This configuration is very slightly shorter than the piecewise-circular
version in~\cite{CKS2}.  As in that version, the core
curve of this tube has discontinuous curvature, for instance at the
``jump point'' where the curve switches from convex to concave.}
\label{fig:brender}
\end{figure}

We now construct a proposed minimizer for a more interesting
Gehring problem---the Borromean rings (see \figr{brender}).
Among the three prime six-crossing links of three components,
the Borromean rings form the one
which is \emph{Brunnian\/}, meaning that if any one component is removed
the remaining components are unlinked.  Milnor's $\mu$--invariant
classifies three-component Brunnian link-homotopy types, and the Borromean
rings are the first nontrivial example.

\begin{figure}[ht!]
\begin{center}
\includegraphics[width=.47\textwidth]{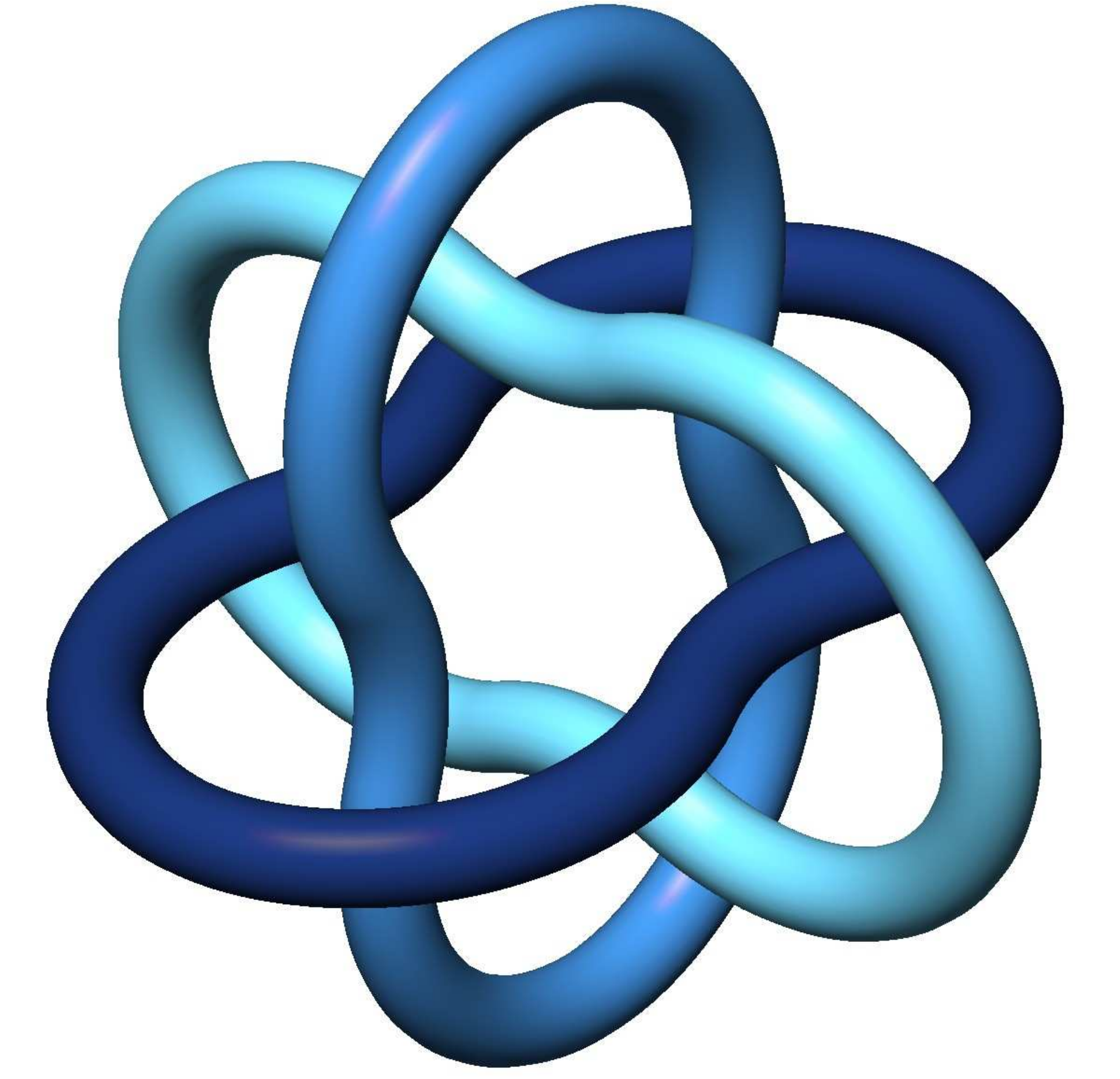}
\includegraphics[width=.47\textwidth]{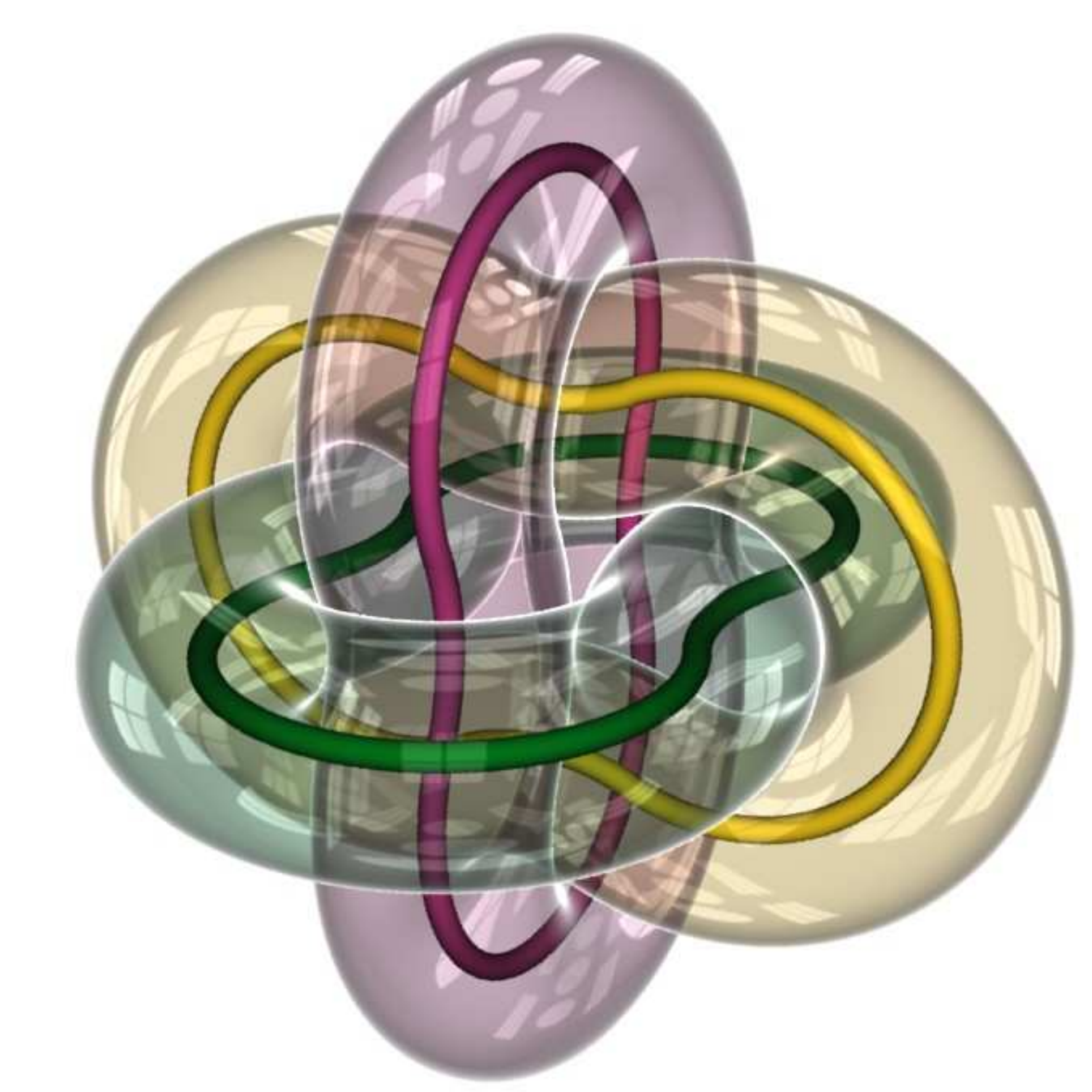}
\end{center}
\caption[The Borromean Rings]
{Two further renderings of the critical configuration~$B_0$
for the Borromean rings reveal more of the structure.
The image on the left, showing thin tubes of diameter~$0.315$,
viewed along an axis of threefold symmetry,
has been adopted as the logo of the International Mathematical Union.
On the right, in a still from the video~\cite{IMU-logo},
we see even thinner tubes inside transparent thick tubes.}
\label{fig:borro-logo}
\end{figure}

In this section, we describe (\thm{borro}) a critical configuration~$B_0$
of the Borromean rings,
shown in \fullref{fig:brender} and \fullref{fig:borro-logo}.
Numerical simulations with Brakke's \emph{Evolver\/}~\cite{evolver}
suggest that this configuration~$B_0$ is in fact the ropelength minimizer
for the Borromean rings.  We will see below that the curvature of~$B_0$ stays
below~$1.534$; this means (as we show in~\cite{CFKSW2}) that~$B_0$ is
also a critical point for length when constrained by the ordinary thickness
measure of~\cite{CKS2} instead of by \gehring thickness.
In~\cite{CKS2}, we described a similar configuration~$B_2$ of the Borromean
rings, built entirely from arcs of unit circles.  \thm{gbalance} shows
that~$B_2$ is not critical, and we compute that~$B_0$ is $0.08 \%$ shorter.

\subsection{Symmetry and convexity}
Our configurations~$B_0$ and~$B_2$ of the Borromean rings are quite similar,
and in particular have the same symmetry and convexity properties, which we
now define.  The three congruent components lie (respectively) in
the three coordinate planes; reflection across any one of these planes
is a symmetry of the link preserving each component.  A further symmetry,
which cyclically permutes the three components, is given by $120^\circ$ rotation
about the $(1,1,1)$ axis; we write this rotation as
$$p\mapsto\tild p\mapsto \ttild p \mapsto p.$$
These symmetries generate the \emph{pyritohedral\/} point group of order~$24$
in~$O(3)$ whose Conway--Thurston orbifold notation \cite{conway3d,conway48}
is~\nstartwo{3}.  Algebraically it is isomorphic to~$\pm A_4$, and in cartesian
coordinates it is most naturally seen as the wreath product $\{\pm1\}\wr C_3$.

\begin{figure}[ht!]\centering
\begin{overpic}[width=4in]{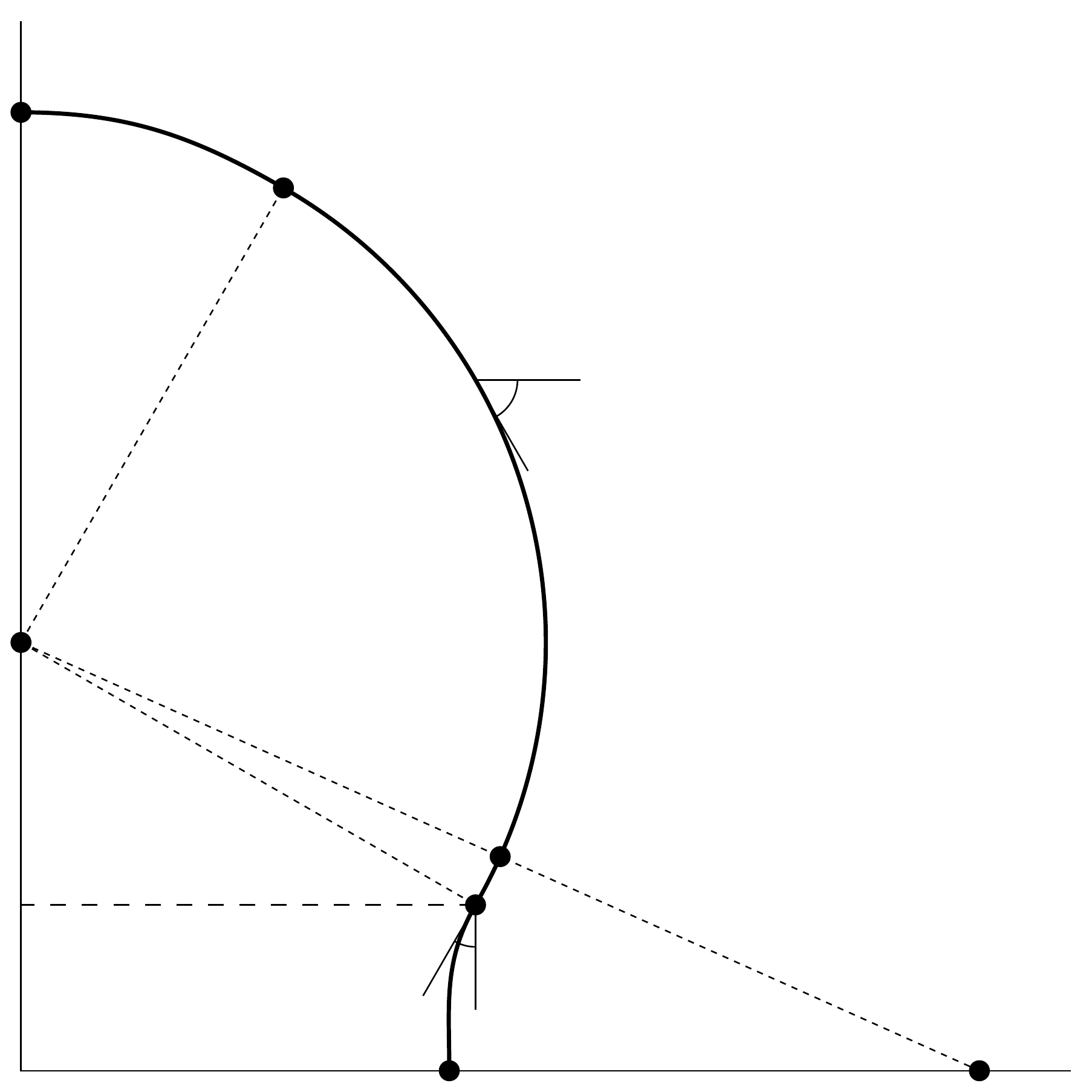}
\small
\put(38,3){$I$}
\put(37,12){$\psi$}
\put(46,15){$J$}
\put(-2,40){$\tild I$}
\put(3,26){$\sigma$}
\put(49,23){$M=\big(\sqrt{1-\rho^2},\rho\big)$}
\put(50,61){$\phi=\arcsin u$}
\put(26,85){$R$}
\put(-3,89){$T$}
\put(89,4){$\ttild T$}
\put(96,3){$x$}
\put(-2,96){$y$}
\end{overpic}
\caption[The generating arc for the Borromean rings]
{Any configuration~$B$ of the Borromean rings with \nstartwo{3}
symmetry is generated by a planar arc~$\arc{IT}$.  We consider
arcs where~$\arc{IJ}$ is concave and~$\arc{JT}$ is convex.
The other points of~$B$ in this quadrant are the rotation
images~$\tild I$ and~$\ttild T$ of~$I$ and~$T$.
In our configurations, there are points~$M$ and~$R$ such
that~$\arc{JR}$ is part of the unit circle around~$\tild I$,
and~$M$ is the midpoint between~$\tild I$ and~$\ttild T$.
The four dotted lines are thus struts of length~$1$.
The height difference from~$J$ to~$\tild I$ is $\sigma=\sin\psi(J)$
as delineated by the horizontal dashed line,
and the coordinates of~$M$ are given
in terms of $\rho=\sin\psi(M)=-\cos\phi(M)$.
}
\label{fig:borroIT}
\end{figure}

Any symmetric configuration of the Borromean rings is the image under
the pyritohedral group of a single embedded arc in the closed
positive quadrant of the $xy$--plane, extending from a point~$I$ on
the $x$--axis to a point~$T$ on the $y$--axis, as shown in \mbox{\figr{borroIT}}.
Conversely, given any such arc~$IT$, its images under \nstartwo{3}
will form a link isotopic to the Borromean rings, as long as~$T$ and~$I$ are
not at the same distance from the origin.  We will assume that
$|I|<|T|$ and will call~$I$ the intip while~$T$ is the tip.
To make the link~$C^1$, the arc~$IT$ must be~$C^1$ and must
meet the axes perpendicularly at its endpoints.

The only other points of the link in this quadrant of the $xy$--plane
are~$\tild I$ and~$\ttild T$; they will be important in the following
discussion.

The arcs~$IT$ of interest to us consist of a small concave piece near the intip
joined to a large convex piece ending at the tip.
That is, there is a jump point $J\in\arc{IT}$ such that the arc~$\arc{IJ}$
is strictly concave, while~$\arc{JT}$ is strictly convex.
As in our discussion of the clasp, we will parametrize~$\arc{IJ}$ by
the angle~$\psi$ (less than~$\tfrac\pi2$) that its tangent vector makes
to the right of the vertical, or by $v=\sin\psi$.
Here~$v$ ranges from~$0$ at~$I$ to some value~$\sigma$ at~$J$,
which will be one of the fundamental parameters for the curves we describe.

Along the convex arc~$\arc{JT}$ we can still define $v=\sin\psi$,
which now decreases from~$\sigma$ through~$0$ to~$-1$.
But we also use the angle $\phi=\tfrac\pi2+\psi$, the angle
above the horizontal made by the tangent vector to~$\arc{JT}$.
Since our curve is~$C^1$, we have $\phi(T)=0$ and
$\phi(J)=\smash{\tfrac\pi2}+\arcsin\sigma$.
In the curves we describe, some initial subarc~$\arc{JR}$ of~$\arc{JT}$ is
part of the unit circle around~$\smash{\tild I}$; we have $\phi(R)\le\smash{\tfrac\pi2}$ so
that along~$\arc{RT}$ we can also use the parameter $u=\sin\phi$.

Finally, to achieve a force balance we will find it necessary
that some point~$M$ along the circle~$\arc{JR}$ has a strut to~$\ttild T$
as well as to~$\tild I$.  This lets us transmit some force
from the large convex arc of one component to the smaller concave arc
of another, indirectly through the third component.
In the $xy$--plane, we find
that~$M$ is the midpoint of the segment~$\tild I \ttild T$,
and thus if we set $\rho := \sin\psi(M)\le \sigma$ we have
$$\tild I=\big(2\rho,0\big), \quad
M=\big(\rho,\tsty\sqrt{1-\rho^2}\big), \quad
\ttild T=\big(0,2\sqrt{1-\rho^2}\big).$$

\subsection{The configuration built from circular arcs}
The configuration~$B_2$ we described in~\cite{CKS2} is
generated by an arc~$\arc{IT}$ of this form.  In~$\smash{B_2}$, we have $R=T$, so that
the entire convex arc~$\arc{JT}$ is part of the unit circle around~$\smash{\tild I}$.
Furthermore, the concave arc~$\arc{IJ}$ is also part of a unit circle,
centered at~$\smash{\ttild T}$. This implies that $M=J$ and $\sigma=\rho=:\rho_2$.
The value~$\rho_2$ is determined by the fact that~$I$ and~$\ttild T$ are
at unit distance, meaning $2\rho_2+1=2\smash{\tsty\sqrt{1-\smash[b]{\rho_2^2}}}$.
As we computed in~\cite{CKS2}, the total length of~$B_2$ is
then $6\pi+24\arcsin\rho_2\approx29.0263$.

This configuration is not balanced (and thus not critical) for \gehring
ropelength.  To balance the curvature forces of the circular arcs, the
fans of struts to their centers would have to carry force proportional
to arclength.  But these struts would then concentrate outward force on
the tips and inward force on the intips; there are no further struts
to balance these forces.  This is like the picture for the naive
clasp---all the force is concentrated on the tips.  As for the clasp,
the tips in the critical configuration will be further apart.

In~\cite{CFKSW2}, we introduce a family of thickness measures with
variable stiffness.  For stiffness~$2$ (meaning that the curves cannot
have osculating circles of diameter less than~$2$) we will see
that~$B_2$ is balanced and hence critical for ropelength.  Because
the circular arcs have exactly the maximum allowed curvature, we will
see that their curvature force need not be balanced pointwise, but only
in total.  Outward strut force on their midpoints (the tips and intips)
can in a sense be spread out to balance the curvature all along the arc.
Because $\rho_2\ne45^\circ$, however, there is an imbalance of
total curvature forces between the convex and concave arcs.
Thus our balancing measure will need an atom of force on the
special colinear struts $\{\tild I, M\}$ and $\{M,\ttild T\}$;
this transmits force from~$\ttild T$ through~$M$ to~$\tild I$.

\subsection{Configurations involving clasp arcs}
To get a balanced configuration~$B_0$ of the Borromean rings, we have to
replace the concave circular arc~$\arc{IJ}$ (and part of the convex arc)
by a tight clasp arc.  Suppose~$\arc{IJ}$ is part of a $\tau$--clasp for
some $\tau\ge\sigma$.  We will now describe a configuration determined by
certain values of our three parameters $$0\le\rho\le\sigma\le\tau\le1,$$
a particular curve of the class illustrated in \figr{borroIT}. 

\begin{figure}[ht!]\centering
\begin{overpic}[scale=1.25]{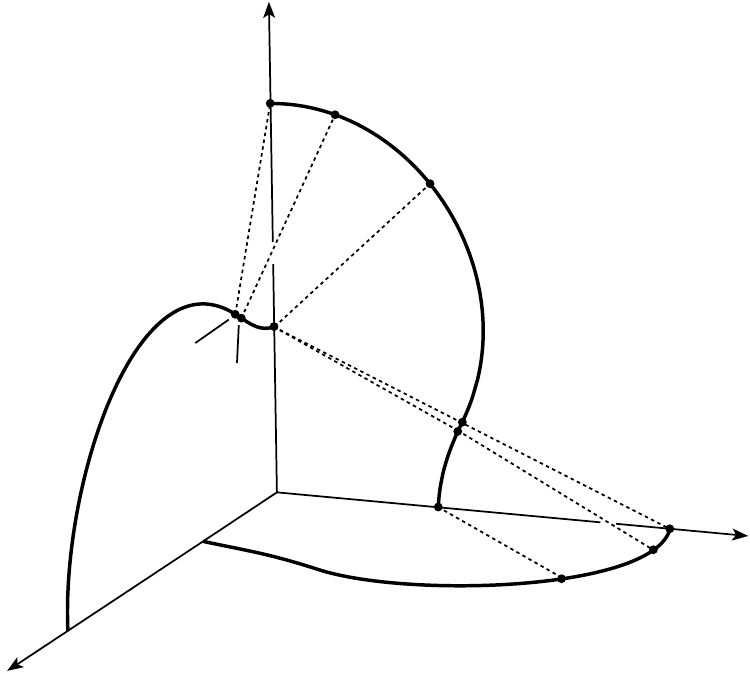}
\small
\put(32,76){$T$}
\put(45.5,75.5){$S$}
\put(58,66){$R$}
\put(63.5,34){$M$}
\put(57,29){$J$}
\put(56,18){$I$}
\put(96,15){$x$}
\put(37.5,84){$y$}
\put(4,5){$z$}
\put(22,40){$\tild M$}
\put(30,37){$\tild J$}
\put(34,41){$\tild I$}
\put(89,21){$\ttild{T}$}
\put(87,13){$\ttild{S}$}
\put(75,8){$\ttild{R}$}
\end{overpic}
\caption[Our Borromean Rings]
{One octant of the critical Borromean rings~$B_0$ consists of three rotated
images of an arc~$\arc{IJMRT}$ of the type shown in \figr{borroIT}.
The dotted lines are struts of length~$1$ connecting the labeled points.
We now describe all other struts to~$\arc{IT}$ in this octant.
Of course, all along the circular arc~$\arc{JMR}$ there are struts to its
center~$\tild{I}$.  Also, between the marked struts are several
one-parameter families of struts, joining two arcs.
The first family joins the conjugate clasp arcs~$\arc{RS}$
and~$\arc{\tild{I}\tild{J}}$; a second family connects~$\arc{ST}$
to the circular arc~$\arc{\tild{J}\tild{M}}$.
The other families are rotated images of these,
connecting~$\arc{JM}$ to~$\arc{\ttild{S}\ttild{T}}$,
and~$\arc{IJ}$ to~$\arc{\ttild{R}\ttild{S}}$.
The struts $\{\ttild{T},M\}$ and $\{M,\tild{I}\}$ are colinear.
To balance~$\arc{IT}$, it is important to
consider also the mirror-image struts across the $xy$--plane.
This figure is an accurate drawing of~$B_0$,
except that we have exaggerated the separation between~$M$ and~$J$: their actual
distance is smaller than the width of the lines used in the picture.
}
\label{fig:bdiagram}
\end{figure}

First, the arc~$\arc{IJ}$ is the piece $v\in[0,\sigma]$ of the
$\tau$--clasp, translated out along the $x$--axis until its tip~$I$ is
at $(2\rho,0,0)$.  Next, $\arc{JMR}$ is an arc of the unit circle
around~$\tild I$, with $v(J)=\sigma$, $v(M)=\rho$ and $u(R)=\tau$.
Note that to get these arcs to match up at~$J$, we will need two conditions on
our parameters~$\rho$, $\sigma$ and~$\tau$.
Finally to define the remaining arc~$\arc{RST}$,
consider the image~$\arc{\tild I\tild J\tild M}$
of~$\arc{IJM}$, rotated into the $yz$--plane.
Then~$\arc{RST}$ is conjugate to~$\arc{\tild I\tild J\tild M}$
in the sense of \lem{perp-planes}:
it is the intersection of the unit-radius tube around
$\arc{\tild I\tild J\tild M}$ with the $xy$--plane,
with~$S$ defined to be the point conjugate to~$\tild J$.
\fullref{fig:bdiagram} shows the arc~$\arc{IT}$ and
its two rotated images, that is,
the part of~$B_0$ lying in the nonnegative orthant in space.

\begin{lemma}\label{lem:bor-match}
For any fixed~$\tau$,
suppose the parameters $0\le\rho\le\sigma\le\tau$ satisfy the two equations
\begin{align}\label{eq:bor-ht}
0 &= 2\rho - \sqrt{1-\sigma^2} +
  \int_{u=0}^{u=\sigma}\!\!\frac{u\du}{\kappa_\tau(u) \sqrt{1-u^2}},\\
\label{eq:bor-c1}
0 &= 1 - (2\rho-\sigma)^2 -
   \frac{1-\sigma^2}{1-\sigma^2(\tau-\sigma)^2},
\end{align}
where~$\kappa_\tau$ is the curvature of the clasp from \thm{clasp}.
Then there is a~$C^1$ and piecewise analytic arc~$\arc{IJMRST}$
as described in the last paragraph.  Its images under the symmetry
group \nstartwo{3} form a configuration $B(\rho,\sigma,\tau)$
of the Borromean rings with \gehring thickness $\gThi=1$.
\end{lemma}

\begin{proof}
As a point on the unit circle~$\arc{JR}$ around~$\tild I$,
the jump point~$J$ has coordinates
$$\big(\sqrt{1-\sigma^2},\, 2\rho-\sigma,\, 0\big).$$
As a point on the $\tau$--clasp $\arc{IJ}$, its coordinates are
$$\Big(2\rho+\int_0^\sigma\!\!\frac{u\du}{\kappa_\tau(u) \sqrt{1-u^2}},\,
  x_\tau(\sigma),\,0\Big).$$
Equating these, using
$$x^2_\tau(\sigma) = 1 - \frac{(1-\sigma^2)}{1-\sigma^2(\tau-\sigma)^2}$$
from \thm{clasp}, gives \eqn{bor-ht} and \eqn{bor-c1}.

If these equations are satisfied, then the position of~$J$ is well-defined,
and~$\arc{IJR}$ is a~$C^1$ arc, meeting the $x$--axis perpendicularly.
The arc~$\arc{RST}$ is the conjugate of~$\arc{IJM}$ and thus
is~$C^1$ by \lem{perp-planes}.  At~$T$, the same lemma shows
it meets the $y$--axis perpendicularly.  At~$R$, the $u=\tau$ base
of the $\tau$--clasp agrees even to second order with the unit circle.

In this configuration, all the struts shown in \figr{bdiagram}
have length~$1$.  If the \gehring thickness were less than~$1$, there would
need to be some shorter strut in this positive octant.  But that
strut would be governed by \lem{perp-planes}, and (rotating to
assume one endpoint is on~$\arc{IT}$) its projection to the
$xy$--plane would be normal to the arc~$\arc{IT}$; the figure
makes it clear that no such strut exists.
\end{proof}

\subsection{The balanced configuration}
Finally, we wish to find the third condition on our parameters~$\rho$,
$\sigma$ and~$\tau$, which will ensure that $B(\rho,\sigma,\tau)$
satisfies the balance criterion of \cor{gehr-bal}.

For most of the struts, it is immediately clear what stress
they need to have in a balancing measure~$\mu$:  The struts
from~$\arc{IJ}$ to~$\arc{\ttild R \ttild S}$, and those
from~$\arc{MR}$ to~$\tild I$ and from~$\arc{RT}$
to~$\arc{\tild I \tild M}$ must be stressed exactly enough
to balance the curvature force of~$\arc{IJ}$ and of~$\arc{MRT}$.
The conjugate clasp arcs~$\arc{IJ}$ and~$\arc{RS}$ exactly
balance each other's curvature forces in this way.

The situation along the short circular arc~$\arc{JM}$ is more complicated.
The struts inward to~$\tild I$ need to balance not only the curvature force
of~$\arc{JM}$ itself, but also the force acting inward on~$\arc{JM}$
from the struts from~$\arc{\ttild S \ttild T}$.  Remember that the measure
needed on these last struts is determined by the curvature of
$\arc{\ttild S\ttild T}$; this in turn determines the measure
needed on the inward struts from~$\arc{JM}$.  We will write this
down explicitly below.  The final condition on our parameters
then comes from a balance of forces at~$\tild I$, where a whole family
of struts converges.

Note that this configuration~$B_0$ of the Borromean rings is the first known
example of a ropelength-critical configuration
in which this sort of transmitted force appears.
Struts impinge on the arc~$\arc{JM}$ from the direction opposite its
own curvature, and transmit their force through that arc.
Without this force transmitted through the (very) short arc~$\arc{JM}$,
the relatively long convex piece~$\arc{\ttild R\ttild T}$ would exert
too much inward force on the relatively short concave piece~$\arc{IJ}$.
Instead, some of this inward force, when transmitted through~$\arc{JM}$,
becomes force \emph{outward\/} on the concave piece~$\arc{\tild I \tild J}$. 
This transmitted force plays the same role in balancing~$B_0$
that the atomic force from~$\ttild T$ through~$M$ to~$\tild I$
played in balancing~$B_2$ for the stiff problem.  But here our
strut measure is absolutely continuous, with no atoms.

To write down the final balancing condition at~$\tild I$,
we begin with an application of \lem{vert-force}:
the total curvature force of~$\arc{JMR}$ and its mirror image
across the $yz$--plane acts on~$\tild I$ downward in the $y$--direction,
with magnitude $$2\big(u(J)-u(R)\big)=2\big(\sqrt{1-\sigma^2}-\tau\big).$$
But the struts from~$\arc{JM}$ carry extra transmitted force.
To determine this, consider the curvature force of an infinitesimal
arc of~$\arc{\ttild{S}\ttild{T}}$ and its mirror image across
the $xy$--plane.  Parametrizing them as usual by~$u$, \lem{vert-force}
tells us the net force, exerted in the negative~$x$ direction,
is~$2\d u$.  This horizontal force is exerted on an infinitesimal
piece of~$\arc{JM}$ and its mirror image across the $xz$--plane.
If we parametrize~$\arc{JM}$ by $v=\sin\psi$, then remembering
that the force on this arc acts perpendicular to the arc, we
see that if its horizontal component is~$\d u$, then its
vertical component is $v\du/\sqrt{1-v^2}$.  This force gets transmitted
through to~$\tild I$.  Because of the symmetry across the $yz$--plane,
of course only the vertical component matters in the end. But this 
symmetry also doubles that vertical force. (\emph{Four\/} copies of the 
arc~$\arc{\ttild{S}\ttild{T}}$ act on~$\tild I$: the original, and reflections
across the $xy$- and $yz$--planes.)  The resultant total transmitted force 
on~$\tild I$ is upward with magnitude
$$ 2\int_{u=0}^{\tau-\sigma}\frac{2v}{\sqrt{1-v^2}}\du.$$
Here the upper limit of integration is $u(S)=\tau-v(J)$
because~$J$ and~$S$ are conjugate points on the $\tau$--clasp.
To make this integral explicit, we need to give the relation
between~$u$ and~$v$; this comes from \lem{perp-planes}.
Along~$\arc{JM}$ we have $y$--coordinate $2\rho-v$,
so the lemma gives
$$u^2 = u(v)^2 := \frac{1-(2\rho-v)^2/v^2}{1-(2\rho-v)^2}.$$
If one wanted, this could be solved to give~$v$ as
the root of a quartic equation in~$\rho$ and~$u$.
Note that $u=0$ at $v=\rho$, as we expect for~$T$ and~$M$.
Plugging in $u=\tau-\sigma$ and $v=\sigma$ (at~$S$ and~$J$)
reproduces \eqn{bor-c1}.

Summarizing, we can write the force-balancing condition at~$\tild I$ as
\begin{equation}\label{eq:bor-bal}
0 = \tau - \sqrt{1-\sigma^2}
 + \int_{v=\rho}^{v=\sigma}\frac{2v}{\sqrt{1-v^2}}\frac{\d u(v)}{\d v}\,\d v
\end{equation}
and so we have proved the following theorem.
\begin{theorem} \label{thm:borro}
Suppose $\rho=\rho_0$, $\sigma=\sigma_0$ and $\tau=\tau_0$
satisfy the three equations \eqn{bor-ht}, \eqn{bor-c1} and \eqn{bor-bal}.
Then the configuration $B_0=B(\rho_0,\sigma_0,\tau_0)$ of the Borromean
rings, constructed as in \lem{bor-match},
is strongly critical for \gehring ropelength.
\end{theorem}

It is easy to solve \eqn{bor-ht} for~$\rho$, or \eqn{bor-c1}
for~$\rho$ or~$\tau$, or \eqn{bor-bal} for~$\tau$, thereby eliminating
one of our three variables.  Then we are
left with two nonlinear integral equations in the other two variables.
While we have not proved formally that a solution to this system
exists, we have solved it numerically to high precision,
both in Mathematica and using \textsc{minpack}~\cite{mp} and \textsc{quadpack}~\cite{qp}.
We obtain
$$ \rho_0 \approx 0.4074218,\quad
 \sigma_0 \approx 0.4177486,\quad
   \tau_0 \approx 0.7561107,$$
where again we follow the standard convention that the error is less
than~$\pm 1$ in the last digit shown.
There is nothing delicate about this solution, since our expressions
vanish to first order at this point.
Numerically it is also clear that this solution is unique.

Using these constants, we compute the length of our critical
Borromean rings~$B_0$ as~$29.0030$.
By comparison, the length of the piecewise circular
Borromean rings~$B_2$ was~$29.0263$.
Thus our critical configuration~$B_0$ beats the naive
circular configuration~$B_2$ by slightly less than
one-tenth of one percent. For comparison, the best lower 
bound known so far~\cite{CKS2}
for the length of the Borromean rings is~$6\pi$.

\begin{figure}[ht!]
\centering
\begin{overpic}[width=9cm]{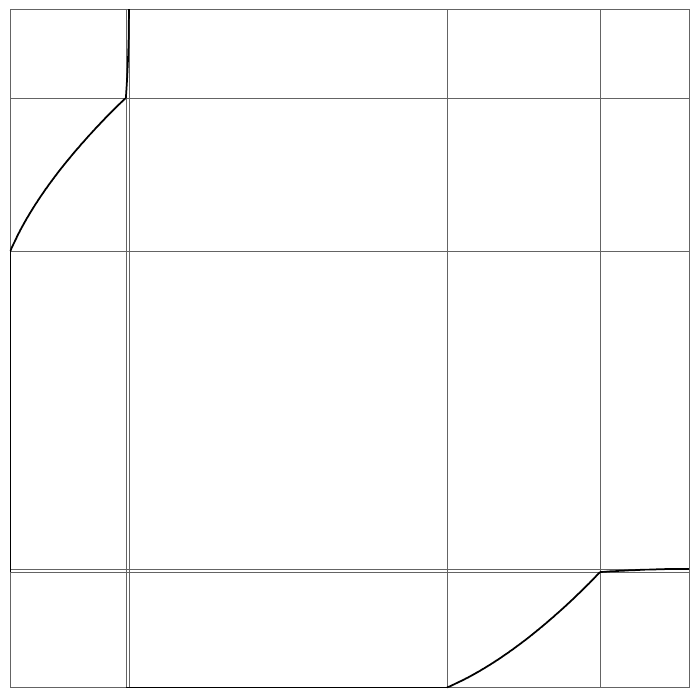}
\small
\put(97,-4){$T$}
\put(84,-4){$S$}
\put(62,-4){$R$}
\put(17.5,-4){$M$}
\put(14,-4){$J$}
\put(0,-4){$I$}
\put(-3,96){$\ttild{T}$}
\put(-3,84.5){$\ttild{S}$}
\put(-3.5,62.5){$\ttild{R}$}
\put(-4.2,19){$\ttild{M}$}
\put(-3.2,13){$\ttild{J}$}
\put(100,19){$\tild{M}$}
\put(100,13){$\tild{J}$}
\put(100,1){$\tild{I}$}
\end{overpic}
\caption[Strut set of the Borromean Rings]
{This picture shows a portion of the strut set of our Borromean
rings, plotted as pairs of arclength coordinates along components of
the link. The horizontal axis represents arclength
along one quadrant of the horizontal component, from~$I$ to~$T$. On
the vertical axis, we plot arclength along quadrants of the other two
components simultaneously.
(This plot accurately depicts the small arclength between~$M$ and~$J$,
in contrast to \figr{bdiagram} where this distance is exaggerated.)
The horizontal segment at the bottom shows the struts from
the circular arc~$\arc{JMR}$ to~$\tild I$; it joins to arcs representing
the families of struts from~$\arc{RS}$ to~$\arc{\tild{I}\tild{J}}$
and from~$\arc{ST}$ to~$\arc{\tild{J}\tild{M}}$.
Symmetrically, the struts to the third component are shown
at the upper left: a vertical segment for the
circle~$\arc{\ttild{J}\ttild{R}}$ around~$I$, and arcs for
the struts from~$\arc{IJM}$ to~$\arc{\ttild{R}\ttild{S}\ttild{T}}$.
Remembering that this square plot should be reflected across
all of its sides to show the complete strut set, we can easily
read off the number of struts coming in
to any point on the curve: two along~$\arc{IJ}$ and~$\arc{RST}$,
three along~$\arc{JM}$ and one along~$\arc{MR}$.
}
\label{fig:bstruts}
\end{figure}

\begin{figure}[ht!]\centering
\begin{overpic}[width=5in]{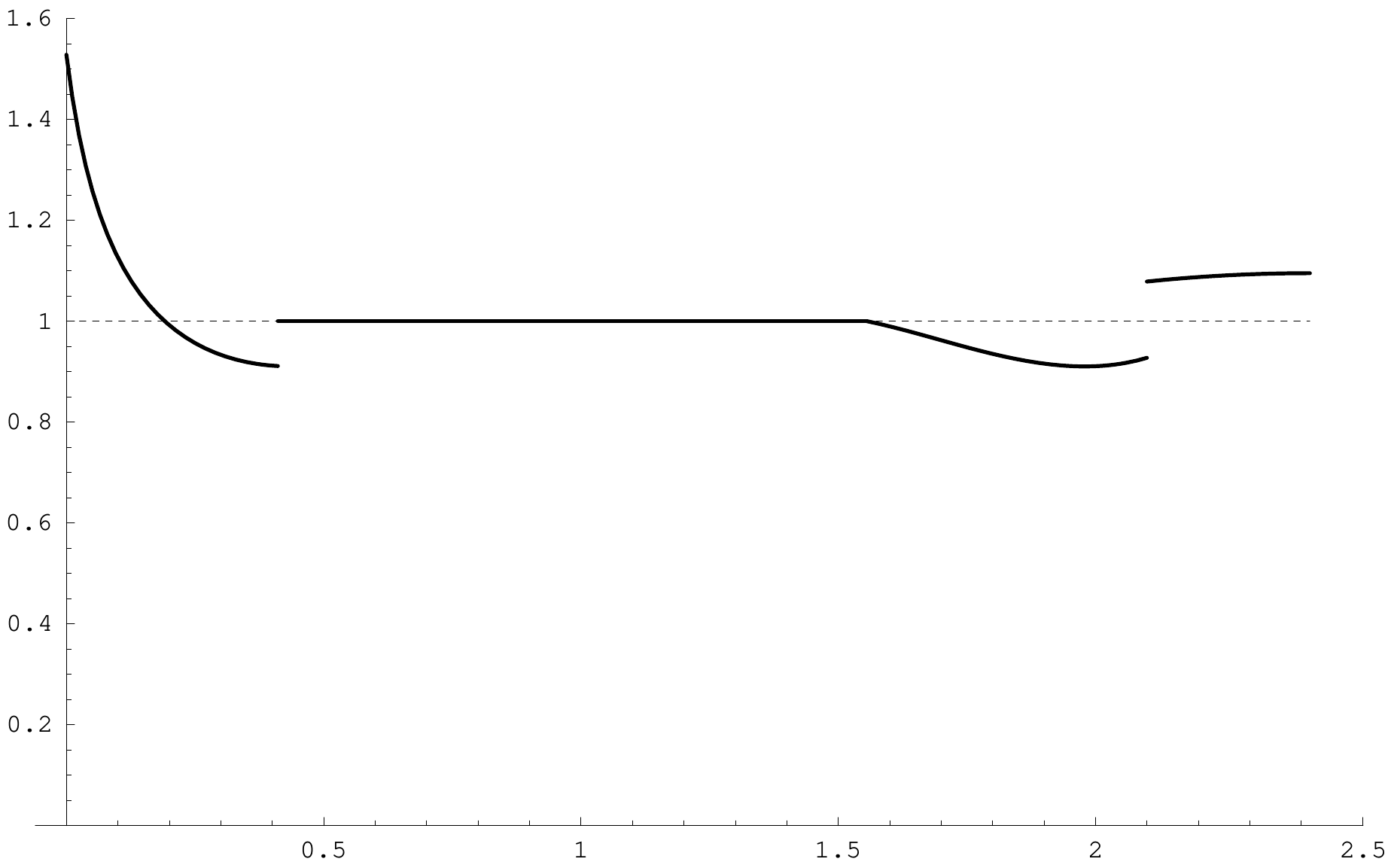}
\small
\put(-2,36){$\kappa$}
\put(50,-2){$s$}
\end{overpic}
\caption[Curvature of the critical Borromean rings]
{The curvature~$\kappa$ of the Borromean rings~$B_0$, plotted
as a function of arclength~$s$ along one quarter of one component
of this critical configuration.  The curvature has its maximum (about~$1.533$)
at the intip~$I$, at $s=0$ in this plot, and then smoothly drops off
to below~$1$.  (This first part could have been plotted negatively,
since this is the concave piece of~$B_0$, but we have chosen to show
the unsigned~$\kappa$ of a space curve.)
After a jump at~$J$, we have $\kappa\equiv1$ along the circular arc~$\arc{JMR}$.
Along the clasp arc~$\arc{RS}$, the curvature drops smoothly from~$1$ and
then rises slightly again, before jumping up above~$1$ at~$S$ and then
increasing to a local maximum at~$T$.}
\label{fig:borrok}
\end{figure}

\fullref{fig:bstruts} shows an arclength plot of the struts
in the Borromean rings.
In \figr{borrok} we plot the curvature of the critical Borromean 
rings~$B_0$ as a function of arclength.  Note that it is discontinuous
only at~$J$ and~$S$.  Each component in~$B_0$ is built of~$14$ analytic
pieces, joined in a $C^{1,1}$ fashion at the symmetric images of the
points~$I$, $J$, $R$ and~$S$.
The maximum curvature (at the intips~$I$)
is $(1 - \tau_0^2)^{-\nicefrac{1}{2}} \approx 1.528$. Therefore~$B_0$ is
also ropelength critical for the standard
ropelength functional of \cite{CKS2}, as we will show in \cite{CFKSW2}.
It is also critical for all the stiff ropelength functionals where the
lower bound~$\lambda$ on the diameter of curvature is less than
$2 \smash{\sqrt{\vphantom{S^+}1-\smash{\tau_0^2}}}\approx 1.3$.

We note that Starostin has described~\cite{starostin-forma} a
configuration~$B_S$
of the Borromean rings with ropelength intermediate between that of
our~$B_2$ and~$B_0$; his configuration replaces the arcs~$\arc{IJ}$
and~$\arc{RT}$ of~$B_2$ by clasp arcs, but does not incorporate the other
features of~$B_0$.  While~$B_S$ can be balanced almost everywhere and
Starostin appears to assume that it is a critical configuration,
in fact it is not balanced at the intips
since it does not satisfy the equivalent of \eqn{bor-bal}.
Thus by \cor{gehr-bal}, $B_S$ is not critical.

\section{Open problems and further directions}
Our work in this paper has been motivated by a simple principle:
that the ideas of rigidity theory for finite frameworks of bars and struts can
be extended to handle mechanisms built from continuous curves of constraints
and contacts. In the simple case of links critical for \gehring ropelength,
this method has already yielded some strong results, such as our
$C^1$--regularity theorem, as well as some surprises like the tight clasp
and the critical Borromean rings. Furthermore we expect that these methods
in general, and our Kuhn--Tucker \thm{ourkt} in particular, will prove
to be useful tools, with applications to a number of
outstanding problems in the geometry and topology of curves and surfaces.

We have mentioned our forthcoming extension of these results~\cite{CFKSW2}
to the classical ropelength problem, where the presence of curvature
constraints and self-contacts of the tube around individual components
makes the situation considerably more challenging.
Our theory of generalized links and obstacles should also be applicable to
the study of packing problems for tubes and surfaces, as when thick rope
is packed into a box~\cite{Kus-pack} (a problem of some interest in molecular
biology \cite{mbms,mmtb}), or when the gray matter of the brain is
folded and pressed against the skull.  We should also mention that while we
have only considered minimizing \emph{length\/} in this paper, our framework
should work equally well for other objective functionals, such as a general
theory of elastic rods with self-contact.

A finite-dimensional duality theorem akin to our Kuhn--Tucker theorem
is one key step in the proof of the Unfolding Theorem
of Connelly, Demaine and Rote~\cite{CDR}:
Any embedded, nonconvex planar polygon admits a motion that preserves
all edge lengths and
strictly increases the distance between any two points on the polygon
not already joined by a straight line of polygon edges.

Our theory allows us to
complete part of the proof of the (conjectured) generalization to
smooth plane curves.  Whether our methods can be made strong
enough to overcome the formidable difficulties involved in
proving a smooth unfolding theorem remains to be seen.

There are several specific open questions suggested by our work above.
\begin{question}
What is the regularity of a critical curve for \gehring ropelength?
Such curves are at worst~$C^1$ and
at best~$C^{1,\nicefrac{2}{3}}$. 
\end{question}
While we have demonstrated \emph{critical\/} configurations of the tight
clasps and Borromean rings, we have not attempted to prove that these
configurations are minimal.
\begin{question}
Are our tight clasps and Borromean rings length-minimal
in their link-homotopy types?
\end{question}
The Euclidean-cone methods of~\cite{CKS2} seem to hold out some hope for
reducing the clasp problem to the case where both curves are planar, but
we have not investigated this line of attack.

\subsubsection*{Acknowledgments}
We gratefully acknowledge many helpful conversations with
our colleagues, including Ted Ashton, Bob Connelly, Erik Demaine,
Elizabeth Denne, Oguz Durumeric, Oscar Gonzalez, Xiao-Song Lin, Marvin
Ortel and Heiko von der Mosel.  In addition, we would like to thank
John Maddocks and the Bernoulli Center at EPFL for hosting some of us during
some of the summer of 2003, Eric Rawdon and Michael Piatek for sharing
their TOROS-minimized link data, and Matt Hoffman for helpful
technical advice on the (extraordinary) Electric Image rendering
system. Some of the other figures in the paper were prepared with
Geomview, \texttt{avn} and RenderMan.
This work was partially supported by the NSF through
grants DMS-99-02397 (to Cantarella), DMS-02-04826 (to Cantarella and Fu),
DMS-00-76085 (to Kusner) and DMS-00-71520 (to Sullivan).

\bibliographystyle{gtart}
\bibliography{link}

\end{document}